\documentclass[journal]{IEEEtran}
\usepackage{cite}

\ifCLASSINFOpdf
\else
\fi
\usepackage{graphicx} 
\usepackage{xcolor}
\usepackage{colortbl} 

\usepackage{amsmath}
\usepackage{amssymb}  
\usepackage[noend]{algpseudocode}
\usepackage{accents}
\usepackage{mathtools}
\usepackage{nicefrac}

\usepackage{algorithm}

\usepackage{booktabs}
\usepackage{multirow,multicol}

\usepackage{url}

\usepackage{tikz}
\usetikzlibrary{fit,matrix,calc}
\usepackage{pgfplots}
\pgfplotsset{compat=1.9}
\newlength\fwidth

\makeatletter
\def\@firstoftwo@second#1#2{%
  \def\temp##1.##2\@nil{##2}%
   \temp#1\@nil}
\newcommand\sref[1]{\expandafter\@setref\csname r@#1\endcsname\@firstoftwo@second{#1}}
\makeatother

\newcommand\MM{M\hspace{-.5mm}M}
\newcommand\MN{M\hspace{-.5mm}N}
\newcommand\twoS{S\hspace{-.3mm}S}

\newcommand{\sdots}{\ifmmode\mathinner{\kern-0.1em\ldotp\kern-0.15em\ldotp\kern-0.15em\ldotp\kern-0.1em}\else.\kern-0.8em.\kern-0.8em.\fi}
\newcommand{\sequal}{\hspace{-.3ex}=\hspace{-.3ex}}

\newcommand{\sminus}{\hspace{-.5ex}-\hspace{-.5ex}}
\newcommand{\splus}{\hspace{-.5ex}+\hspace{-.5ex}}
\newcommand{\shortin}{\hspace{-.3ex}\in\hspace{-.3ex}}

\newcommand*\circled[1]{\raisebox{1pt}{\tikz[baseline=(char.base)]{\node[shape=circle,draw,inner sep=1pt] (char) {\tiny #1};}}}

\newcommand{\bmat}[1]{\begin{bmatrix} #1 \end{bmatrix}}

\usepackage{accents}
\makeatletter
\newcommand{\sinternal}{%
  \hbox{\scalebox{.2}{$\blacksquare$}}}
\DeclareRobustCommand{\internal}{\accentset{\sinternal}}
\makeatother

\makeatletter
\newcommand{\sexternal}{%
  \hbox{\scalebox{.25}{$\square$}}}
\DeclareRobustCommand{\external}{\accentset{\sexternal}}
\makeatother

\DeclareMathOperator*{\cost}{\operatorname{c}}
\DeclareMathOperator*{\proj}{\Pi}
\DeclareMathOperator*{\diag}{\operatorname{diag}}

\definecolor{lightblue}{rgb}{0,         0.4470,    0.7410}
\definecolor{lightred}{rgb}{0.8500,    0.3250,    0.0980}
\definecolor{lightgreen}{rgb}{0.4660,    0.6740,    0.1880}
\definecolor{darkpurple}{RGB}{119,63,155}

\newif\ifArxiv
\Arxivtrue

\newcommand{\eqlabel}[1]{\addtocounter{equation}{-1}\refstepcounter{equation}\label{#1}}

\usepackage{amsthm} 
\theoremstyle{plain} 
\newtheorem{proposition}{Proposition}

\newtheorem{definition}{Definition}
\newtheorem{assumption}{Assumption}

\theoremstyle{definition}
\newtheorem{example}{Example}

\usepackage[%
  pdftitle={ADMM for Exploiting Structure in MPC Problems},
  pdfauthor={Felix Rey, Peter Hokayem, John Lygeros},
  pdfcreator={PDFlatex},
  pdfsubject={},
  pdfkeywords={},
  hypertexnames=false, 
  bookmarks=true,
  colorlinks=true,
  linkcolor=black,
  citecolor=black,
  pagebackref=false,
  breaklinks=true,
  raiselinks=true,
  urlcolor=black]{hyperref}

\begin{document}
%
\title{ADMM for Exploiting Structure in MPC Problems}
%
%
%

\author{Felix~Rey,~\IEEEmembership{Member,~IEEE,}
        Peter~Hokayem,~\IEEEmembership{Member,~IEEE,}
        and~John~Lygeros,~\IEEEmembership{Fellow,~IEEE}
\ifArxiv
\else
\thanks{Manuscript received August xx, 2018; revised xxxx xx, 2019.}
\fi
\thanks{This work was supported by ABB Corporate Research under Grant 2017-1224/01. We thank Michael Cantoni from the University of Melbourne for bringing the well-suited application to cascade systems to our attention.}
\thanks{F.~Rey {\tt\small rey*} and J. Lygeros {\tt\small lygeros*} are with the Automatic Control Laboratory at ETH Zurich, Switzerland, {\tt\small *@control.ee.ethz.ch}.}
\thanks{Peter Hokayem {\tt\small peter}{\tt\small.al-hokayem}{\tt\small@ch.abb.com}  is with ABB, Switzerland.}
}
\maketitle

\begin{abstract}
We consider a model predictive control (MPC) setting, where we use the alternating direction method of multipliers (ADMM) to exploit problem structure. We take advantage of interacting components in the controlled system by decomposing its dynamics with virtual subsystems and virtual inputs. We introduce subsystem-individual penalty parameters together with optimal selection techniques. Further,
we propose a novel measure of system structure, which we call separation tendency. For a sufficiently structured system, the resulting structure-exploiting method has the following characteristics: $(i)$~its computational complexity scales favorably with the problem size; $(ii)$ it is highly parallelizable; $(iii)$ it is highly adaptable to the problem at hand;  and $(iv)$, even for a single-thread implementation, it improves the overall performance. 
We show a simulation study for cascade systems and compare the new method to conventional ADMM.
\end{abstract}


%
\IEEEpeerreviewmaketitle

\section{Introduction}
\IEEEPARstart{I}{n} early
 applications, model predictive control (MPC)~{\cite{maciejowski2002predictive,camacho2013model}} was used in large and cost-intensive processes, for example in the chemical industry~\cite{morari1999model}. In such settings, where the expense of the control hardware in the overall process is small, we can use a generously-sized computation device that is capable of hosting a general-purpose solver. In contrast, the rapid advance of information technology brings MPC-based control techniques to mass production~\cite{di2012industry}, e.g., in automotive industries and consumer electronics.  In such large-volume production settings, the pressure on cost-per-unit calls for high efficiency, which we can achieve through tailored hardware and specialized algorithms.
 We focus on the optimization algorithm, and we use the alternating direction method of multipliers (ADMM)~\cite{glowinski1975approximation,boyd2011distributed}, a first-order method for solving convex problems.
On the hardware side, we consider embedded platforms, such as field-programmable gate arrays (FPGAs) or application-specific integrated circuits (ASICs). ADMM suits embedded devices as it performs only simple and numerically stable operations~\cite{jerez2014embedded}. 
Our goal is to increase the synergy between the algorithm, platform, and problem, by adapting ADMM  to fit the controlled system. Hence, instead of pursuing a generalist approach that most off-the-shelve solvers provide, we develop a specialization strategy that improves the performance in a specific MPC setting. 

\begin{figure}[t]
   \centering
  \includegraphics[width=.93\columnwidth]{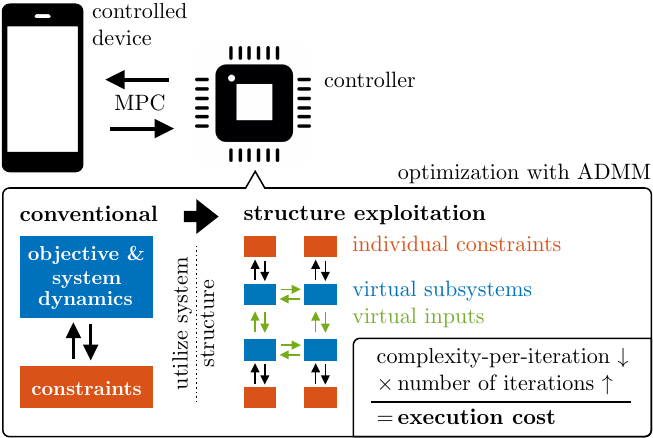}
  \caption{MPC setup with a controller running structure-exploiting ADMM.}
   \label{fig:setup}
 \end{figure}

  We solve MPC problems that are generally composed of a control objective, system dynamics, and additional state and input constraints. When we conventionally apply ADMM~\cite{raghunathan2014optimal}, the algorithm mediates between the system dynamics and the additional constraints. 
  We go beyond this formulation by exploiting structure in the controlled system through a decomposition into \textit{virtual subsystems}, connected through \textit{virtual inputs}. The tailored algorithm then mediates between all subsystems and reassembles the full system dynamics only in convergence.  Fig.~\ref{fig:setup} illustrates the setup, where we also define the overall cost for executing the algorithm. Our numerical results show that the decomposition typically increases the required number of iterations compared to the conventional application of ADMM. However, the structure-exploiting method reduces the complexity of each iteration and has a large parallelization potential. Overall, if the controlled system is sufficiently structured, structure exploitation  reduces the execution cost.\looseness=-1

 

Beside ADMM, Fast Gradient Methods (FGMs)~\cite{jerez2013embedded} are a common choice for embedded optimization~\cite{peyrl2014parallel}. 
ADMM permits more design freedom and is efficiently applicable in a broader problem range, e.g., for state constraints or convex but not strongly convex objective functions~\cite{boyd2011distributed}.
Examples of problem-specific adaptation through ADMM are~\cite{vaya2014decentralized} for electric vehicle charging; \cite{kang2015decomposition} for scenario-based stochastic programs; and \cite{Rey2017d} for coordinated energy reserve bidding. In these publications, a common theme is the problem separation over a coupling constraint. While we employ a similar separation, we utilize a virtual decomposition that is not present in the original problem. In~\cite{deng2017parallel,chen2016direct,mateos2010distributed,chang2015multi}, general distributed ADMM formulations are discussed, which do not focus on MPC and do not exploit structure in the problem data.  As opposed to potentially non-convergent multi-block extensions~\cite{deng2017parallel,chen2016direct}, our method remains part of the standard ADMM family~\cite{boyd2011distributed}, which guarantees convergence and makes our approach amenable to existing results, e.g., infeasibility detection~\cite{stellato2017osqp,raghunathan2014admm},  prescaling~\cite{GhadimiOptimalParameter,Rey2016,stathopoulos2016operator}, and over-relaxation~\cite{giselsson2017linear,boyd2011distributed}.
Techniques that focus on the MPC setup structure (as opposed to the system structure and the ADMM utilization) are also compatible, e.g., move-blocking~\cite{wang2010fast} and warm-starting~\cite{wang2010fast,stathopoulos2016operator}. 

Two ADMM formulations~\cite{raghunathan2014optimal,GhadimiOptimalParameter} are widely used in MPC. For the formulation in~\cite{GhadimiOptimalParameter}, the problem needs to be condensed through constraint elimination, which also eliminates structure. With the method in~\cite{raghunathan2014optimal}, no condensing is required.
The formulations also differ in the addressed problem type. In~\cite{raghunathan2014optimal}, state constraints are possible, and the formulation is particularly suited for easily-projectable constraint sets. We build on this formulation as it preserves the problem structure.

Our method is the result of ongoing research~\cite{Rey2017a,Rey2017b}. In~\cite{Rey2017a}, we use the basic idea of system structure exploitation for optimization problems that arise from controlling power conversion systems. We generalize the approach to an MPC framework in~\cite{Rey2017b}, where we also adopt the name structure-exploiting ADMM. In the present work, we advance the overall framework and we present two major novelties: We introduce subsystem-individual penalty parameters with an optimal selection technique, and we propose a novel  measure of system structure that we call separation tendency.

\vspace{\baselineskip}
\paragraph*{Notation}
We denote dimensions with matching non-italic symbols, e.g., $x \in \mathbb R^{\mathrm x}$. The identity matrix is $I_n\in \mathbb R^{n\times n}$, and $0_{n\times m}$, $1_{n \times m}$ are $n \times m$ matrices with all elements $0$ or $1$. We omit the subscripts if the dimension is clear from context. We write $\|x\|_{Q}^2$ for $x^\top Q x$.
We use $[A_{ij}]$ to concatenate matrices $A_{ij}$ along row~$i$ and column~$j$, as well as vertical $[A_{1j};A_{2j}]$, horizontal $[A_{i1},A_{i2}]$, and diagonal $\operatorname{diag}(A_{11}, A_{22})$ \mbox{concatenation}. We denote a sequence of elements $x_i$ with $\{x_i\}_{i=1,2,\dots}$, where we again may omit the range of $i$.

\section{MPC Formulation and System Structure}
We consider  the MPC problem
\begin{IEEEeqnarray}{rCll}
\IEEEyesnumber\eqlabel{eq:mpcprog} \IEEEyessubnumber*
\min_{\{x^{k+1},u^k\}} &\quad& \sum\nolimits_{k=1}^N \big( \tfrac12 \|x^{k+1} - r_x^k  \|_{Q}^2 \rlap{$  + \tfrac12\|u^{k} - r_u^k \|_{R}^2 \big)$} &  \hspace{3cm} \label{eq:mpcprogobjective}\\
\text{s.t.}&& x^{k+1} = Ax^k + Bu^k  &\forall\,k=1,\sdots, N  \label{eq:initsys}\\
&& (x^{k+1},u^k) \in \mathcal X \times \mathcal U &\forall\,k=1,\sdots, N,  \label{eq:mpcprog_ineq}
\end{IEEEeqnarray}
with prediction horizon $N$, state \mbox{$x^k \shortin R^{\mathrm x}$}, input \mbox{$u^k \shortin R^{\mathrm u}$}, dynamics matrix $A$, input matrix $B$, tracking references~$r^k_x, r^k_u$, symmetric weights $Q,R$, and constraint sets~$\mathcal X, \mathcal U$. 
We use Assumption~\ref{ass:convexity} throughout the rest of the paper.

\begin{assumption}
\label{ass:convexity}
  Problem~\eqref{eq:mpcprog} is convex and feasible. 
\end{assumption}

Problem~\eqref{eq:mpcprog} is convex if $Q$, $R$ are positive semidefinite and $\mathcal X$, $\mathcal U$ are convex. The problem is feasible if there exist trajectories $\{x^{k+1},u^k\}_{k=1,\dots,N}$ that satisfy~\eqref{eq:initsys}, \eqref{eq:mpcprog_ineq}. 
 Our approach is particularly suited for problems where a projection onto~$\mathcal X, \mathcal U$ is computationally cheap, as these projections will be used repeatedly.  

\subsection{State and Input Partition}
\label{subsec:partition}
\begin{figure}[t]
   \centering
  \includegraphics[width=.9\columnwidth]{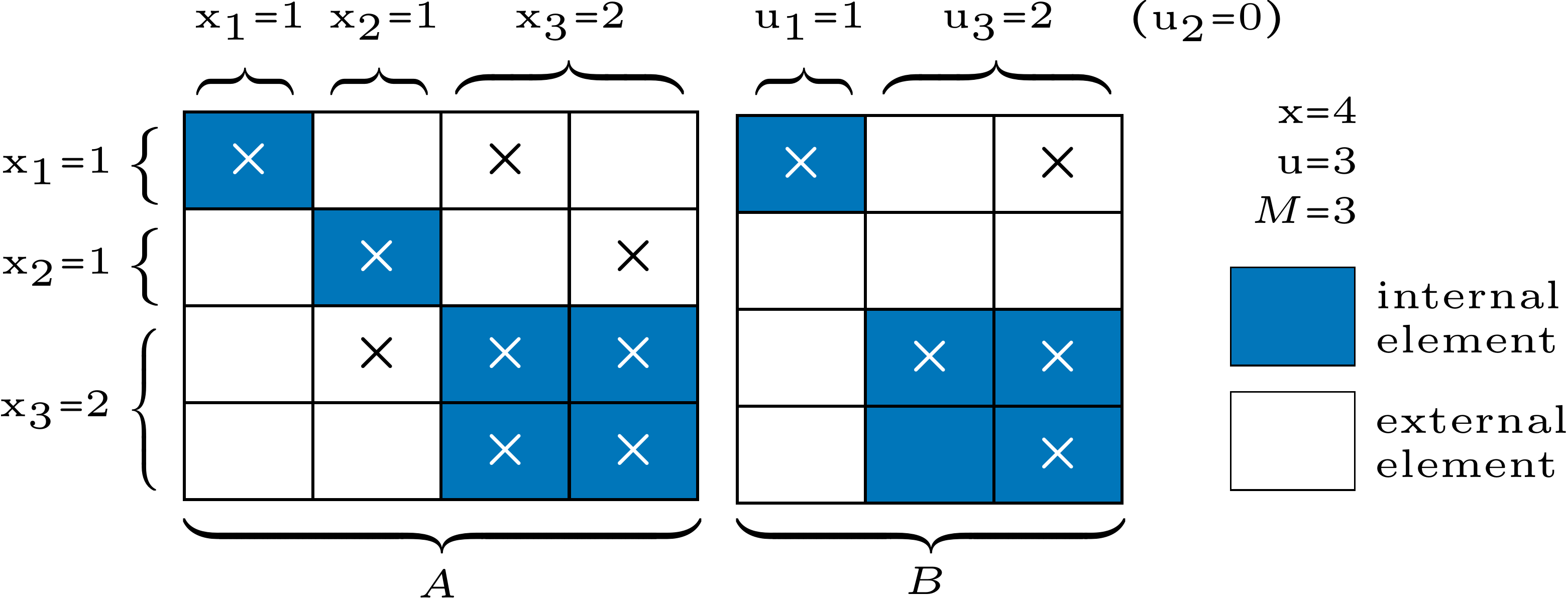}
  \caption{Illustration of a structured system, where $\times$ is a placeholder for any non-zero element. We choose the partition to best-resemble dense components, indicated by the matrix occupation pattern. A suitable partition leaves few non-zero external elements, while it still decomposes the system.}
  \label{fig:partitionExample}
 \end{figure}

 Our method exploits interacting components in the system~$(A,B)$. 
 We formalize the presence of such components with the state and input partition $x^k\sequal[x^k_1;\sdots;x^k_M]$, $u^k\sequal[u^k_1;\sdots;u^k_M]$. We use $x^k_i\in\mathbb R^{\mathrm x_i}$, $u_i^k \in \mathbb R^{\mathrm u_i}$ with \mbox{$\mathrm x_i \in \mathbb N_{>0}$}, \mbox{$\mathrm u_i \in \mathbb N_{\geq 0}$}, $\sum_i \mathrm x_i \sequal \mathrm x$, and $\sum_i \mathrm u_i \sequal \mathrm u$. 
We assume that states and inputs are ordered already, such that only consecutive elements are grouped. Hence, the dimensions $\{\mathrm x_i, \mathrm u_i\}_{i=1,\dots,M}$ define the partition. 
In Fig.~\ref{fig:partitionExample}, we show a system with an exemplary component pattern and introduce the notion of internal and external elements. 
The partition decomposes~$A,B$ into {submatrices} $A_{ij}\sequal \mathbb R^{\mathrm x_i\times\mathrm x_j}$, $B_{ij}\sequal \mathbb R^{\mathrm x_i\times\mathrm u_j}$ such that $A\sequal [A_{ij}]$, $B\sequal [B_{ij}]$. As $\mathrm u_i\geq 0$, empty submatrices can result from the decomposition of~$B$.
Furthermore, we decompose~$A,B$ into a sum of \textit{internal} matrices~$\internal{A}, \internal{B}$ and \textit{external} matrices~$\external{A}, \external{B}$, where $A \sequal \internal{A} \hspace{-0.1ex}+\hspace{-0.1ex} \external{A}$ with $\internal{A} \sequal \diag({A_{11},\sdots,A_{\MM}})$ and $B \sequal \internal{B} \hspace{-0.1ex}+\hspace{-0.1ex} \external{B}$ with $\internal{B} \sequal \diag({B_{11},\sdots,B_{\MM}})$.

Our aim is to partition Problem~\eqref{eq:mpcprog} along the same lines as the system. Towards this, we require the following assumption.
\begin{assumption}
\label{ass.admissible} The partition $\{\mathrm x_i, \mathrm u_i\}$ can be used to decompose~\eqref{eq:mpcprogobjective} and~\eqref{eq:mpcprog_ineq}, i.e., there exist $Q_i$, $R_i$, $\mathcal X_i$, $\mathcal U_i$, $i=1,\sdots,M$ such that~{$Q \sequal \diag(Q_1,\sdots,Q_M)$}, 
    {$R \sequal \diag(R_1,\sdots,R_M)$},
    {$\mathcal X \sequal \{ x^k\, |\, x_i^k\in  \mathcal X_i \,\,\,\forall i\}$}, and
    $\mathcal U \sequal \{ u^k \,|\, u_i^k\in  \mathcal U_i \,\,\,\forall i\}$. 
\end{assumption}

We say the partition is \emph{admissible} for~\eqref{eq:mpcprog} if Assumption~\ref{ass.admissible} is satisfied. 
    Any partition is admissible if~$Q,R$ are diagonal and~$\mathcal X,\mathcal U$ are separable. Conversely, the trivial partition~$M\sequal 1$ is admissible for any problem.
    If the external part of the partition is sparse, our approach will lead to computational benefits. In such a case, we call the respective system \emph{structured}. We present a better-quantified structure measure in Section~\ref{sec:separationTendency}. 
If a suitable partition cannot be found by inspection, we can use spectral clustering methods \cite{ageev2000approximation,frieze1997improved,hespanha2004efficient}, where~\cite{hespanha2004efficient} also provides a suitable ordering for states and inputs.


\subsection{Virtual Inputs and Subsystems}
By using an admissible partition $\{\mathrm x_i, \mathrm u_i\}$, we rewrite system~\eqref{eq:initsys} as
  $x^{k+1} \sequal \internal{A}x^k + \internal{B}u^k + v^k $
with the virtual input 
\begin{IEEEeqnarray}{rCl}
\label{eq:defvirtualinput}
v^k = \external{A}x^k + \external{B}u^k,
\end{IEEEeqnarray}
which represents the external coupling. As $\internal{A}, \internal{B}$ are block-diagonal, we can make the decomposition more explicit by rewriting the system  as a collection of virtual subsystems
\begin{IEEEeqnarray}{rCl}
\label{eq:defvirtualsubsys}
x_i^{k+1} = A_{ii}x_i^k + B_{ii}u_i^k + v_i^k, \quad i=1,\sdots,M,
\end{IEEEeqnarray}
 where~$v^k$ adopts the partition of~$x^k$. If $\external{A}, \external{B}$ are sparse, then the virtual input $v^k_i \in \mathbb R^{\mathrm x_i}$ can take values in a lower dimensional space than suggested by its dimension~$\mathrm x_i$. To make this explicit, we write the components in~\eqref{eq:defvirtualinput} as
  \begin{IEEEeqnarray}{rCl}
\label{eq:defvirtualinput2}
  v^k_i= \sum\nolimits_{j}  A_{ij} x_j^k + \sum\nolimits_{j}  B_{ij} u_j^k,
\end{IEEEeqnarray}
  where $j\shortin \{1,\sdots, M \} \backslash \{i\}$. We pick a matrix \mbox{$W_i \in \mathbb R^{\mathrm x_i \times \mathrm w_i}$} with $\mathrm w_i \leq \mathrm x_i$ such that its columns form a range space basis for the concatenated matrix
    \begin{IEEEeqnarray}{rCl}
    \label{eq:concetenatedMatrix}
  \left[\left\{ A_{ij}, B_{ij} \right\}_{ j\in \{1,\dots, M \} \backslash \{i\} } \right] \in \mathbb R^{\mathrm x_i \times (\mathrm x + \mathrm u - \mathrm x_i - \mathrm u_i)}, 
\end{IEEEeqnarray}
  which we obtain by writing~\eqref{eq:defvirtualinput2} as a single matrix-vector operation. Hence,~\mbox{$\mathrm w_i$} is the row rank of~\eqref{eq:concetenatedMatrix}, and the range space of~$W_i$ contains all values that~$v_i^k$ can attain.  We then introduce the dimension-reduced virtual input \mbox{$w^k_i\in \mathbb R^{\mathrm w_i}$} by replacing~$v_i^k$ with~$W_iw_i^k$.

\subsection{Partitioned Problem}
By using an admissible partition $\{\mathrm x_i, \mathrm u_i\}$,  we write~\eqref{eq:mpcprog} as
\begin{IEEEeqnarray}{rCll}
\IEEEyesnumber\eqlabel{eq:mpcprogpartitioned} \IEEEyessubnumber*
\min_{\{x^{k+1},u^k,w^k\} } &\,& \sum\limits_{i,k=1}^{M,N}\hspace{-0.7ex} \Big( \tfrac12 \|x_i^{k+1} - r^k_{x_i} \|_{Q_i}^2  \hspace{-0.1ex}+ \hspace{-0.1ex}\rlap{$\tfrac12 \|u_i^{k} - r^k_{u_i} \|_{R_i}^2 \Big)$} \\
\text{s.t.} &\,\,\,& x_i^{k+1} =  A_{ii} x_i^k +  B_{ii} u_i^k +  W_i w_i^k \quad&\forall \,(i,k) \quad\quad\,\\
&& (x_i^{k+1},u_i^k) \in \mathcal X_i \times \mathcal U_i \quad \quad&\forall \,(i,k) \label{eq:mpcprogpartitioned_constraints}\\
&&  W w^k = \external{A}x^k + \external{B}u^k  &\forall \,k,  \label{eq:mpcprogpartitioned_coupling}
\end{IEEEeqnarray}
where $w^k\hspace{-.19ex}=\hspace{-.19ex}[w^k_1;\sdots;w^k_{\hspace{-.19ex}M}] \hspace{-.19ex} \in \hspace{-.19ex} \mathbb R^{\mathrm w}$\hspace{-.19ex} and $W \hspace{-.19ex}=\hspace{-.19ex} \diag(W_1,\sdots,W_{\hspace{-.19ex}M})$. Furthermore, with the stacked variables $y_i^{k} \sequal [u^k_i;w^k_i;x^{k+1}_i]$, \mbox{$ y_i \sequal [y_i^1;\sdots;y_i^{N}]\in \mathbb R^{\mathrm y_i}$}, and \mbox{$ y \sequal [ y_1;\sdots; y_M]\in \mathbb R^{\mathrm y}$}, we obtain 
\begin{IEEEeqnarray}{rCllr}
\IEEEyesnumber\eqlabel{eq:mpcprogfinal} \IEEEyessubnumber*
\min\limits_{y} &\quad& \sum\nolimits_{i=1}^M \big( \tfrac12 y_i^\top {}& \mathcal Q_i y_i \rlap{$\,+\, q_i^\top y_i +  K_i \bigr)$} \hspace{5mm}& \text{ (objectives)}\hspace{5mm} \label{mpcfinal_objective} \\[-0.2ex]
\text{s.t.} && C_i y_i=c_i& \forall i  &\text{ (dynamics)}\hspace{5mm} \\
&&  y_i \in \mathcal Y_i & \forall i &\text{ (constraints)}\hspace{5mm}  \label{eq:mpcfinal-subsysBounds} \\
&& D  y  = d.  & \hspace{3cm}&\text{(coupling)}\rlap{}\hspace{5mm}  \label{eq:mpcfinal-coupling} 
\end{IEEEeqnarray}
We list the definitions of $\mathcal Q_i, q_i, K_i, C_i, c_i, \mathcal Y_i, D$, and $ d$ in Appendix~\ref{app:stacking}. Problem~\eqref{eq:mpcprogfinal} is equivalent to Problem~\eqref{eq:mpcprog}. 

\subsection{Conventional ADMM}
\label{sec:conventionalADMM}
Our approach extends the conventional ADMM formulation in~\cite{raghunathan2014optimal}, which considers the unpartitioned problem
\begin{IEEEeqnarray}{rCllr}
\IEEEyesnumber\eqlabel{eq:mpcunpartitioned} \IEEEyessubnumber*
\min\limits_{y} &\quad& \tfrac12 y^\top \mathcal Q y + q^\top  y + K\\
\text{s.t.} && C y=c   \\
&&  y \in \mathcal Y.\label{eq:mpcunpartitioned_ineq}
\end{IEEEeqnarray}
For $M\sequal 1$, Problem~\eqref{eq:mpcunpartitioned} is equivalent to Problem~\eqref{eq:mpcprogfinal}. As in~\cite{raghunathan2014optimal}, we rewrite~\eqref{eq:mpcunpartitioned_ineq} with 
$y \sequal \zeta$, $\zeta  \shortin \mathcal Y$,
where $\zeta$ is a duplicate of~$y$. 
As shown in Algorithm~\ref{alg:conventionalADMM}, ADMM then addresses the parts depending on~$y$ and~$\zeta$ alternatingly.
We use an orthogonal projection $\Pi$ as in~\cite{raghunathan2014optimal,parikh2014proximal}, a scaled Lagrange multiplier $\lambda$ which is  associated to the constraint~$y=\zeta$, and a user-defined penalty parameter~\mbox{$\rho$}  that influences the convergence speed~\cite{boyd2011distributed}.
 It is shown in~\cite{raghunathan2014optimal} that Algorithm~\ref{alg:conventionalADMM} with $\rho\hspace{-.0ex}>\hspace{-.0ex}0$ converges to a fixed point that is optimal for~\eqref{eq:mpcunpartitioned}. Suitable initialization and termination techniques are discussed in~\cite{boyd2011distributed}, and the behavior for infeasible problems is discussed in~\cite{raghunathan2014admm,banjac2017infeasibility}.

\begin{algorithm}
\caption{Conventional ADMM}\label{alg:conventionalADMM}
\begin{algorithmic}[1]
\Statex \hspace{-1em}\textbf{repeat}

\algrenewcommand{\alglinenumber}[1]{\hspace{1.4ex}\footnotesize\color{black}\circled{1.1}}
\State  $y \leftarrow \arg\hspace{-.3ex}\min\limits_{y} \, \tfrac12 {y^\top} \hspace{-.5ex} \mathcal Q {y} + q^\top \hspace{-.5ex} {y}
+ 
\tfrac{\rho}{2} \| {y} - {\zeta} - \lambda \|_2^2
$
\Statex \hspace{2.11cm}  $\text{s.t.} \,\, C {y}=c$
\vspace{.8ex}

\algrenewcommand{\alglinenumber}[1]{\hspace{1.4ex}\footnotesize\color{black}\circled{1.2}}
\State  $ {\zeta} \leftarrow \proj_{\mathcal Y} \left(  {y}   -  \lambda \right)$  
\vspace{0.8ex}

\algrenewcommand{\alglinenumber}[1]{\hspace{1.4ex}\footnotesize\color{black}\circled{1.3}}
\State  $ \lambda \leftarrow   \lambda-  (y-\zeta)$
\end{algorithmic}
\end{algorithm}

\section{Structure-Exploiting ADMM}
\label{sec:structureexplADMM}

The structure-exploiting algorithm utilizes the problem partition and therefore takes advantage of the system structure.

\subsection{Main Algorithm}
To solve Problem~\eqref{eq:mpcprogfinal} with $M>1$, we add (additional to~$\zeta$) the second duplicate~$\epsilon$. Equivalent to~\eqref{eq:mpcprogfinal}, we obtain
\begin{IEEEeqnarray}{rCllrl}
\IEEEyesnumber\eqlabel{eq:mpcprogsplitting} \IEEEyessubnumber*
\min\limits_{\textcolor{lightblue}{y}, \textcolor{lightred}{\zeta}, \textcolor{lightgreen}{\epsilon} } &\quad&   \sum\nolimits_{i=1}^{M} \big( \tfrac12  \textcolor{lightblue}{y_i^{\textcolor{black}{\top}}} & \mathcal Q_i \textcolor{lightblue}{y_i}  \rlap{$ + q_i^\top  \textcolor{lightblue}{y_i} \color{black} \splus K_i \bigr)$} \hspace{1.15cm} &  \color{lightblue}\text{\raisebox{-1.5mm}{{(individual objec-}}}  \\[-0.2ex]
\text{s.t.} &&  C_i \textcolor{lightblue}{y_i}=c_i& \forall i  &\color{lightblue} \text{{tives and dynamics)}} \\
&&  \textcolor{lightred}{\zeta_i} \in \mathcal Y_i & \forall i &\color{lightred} \text{{(individual constraints)}} \label{eq:mpcprogsplitting_individual}\\
&&  D  \textcolor{lightgreen}{\epsilon}  = d  & & \color{lightgreen}\text{{(coupling)}} \label{eq:mpcprogsplitting_coupling}   \\
&& \textcolor{lightblue}{y} = \textcolor{lightred}{\zeta} = \textcolor{lightgreen}{\epsilon}. && \llap{\text{(variable duplication)}} & \quad\quad \label{eq:mpcprogsplitting_consenus}
\end{IEEEeqnarray}
 We show the resulting structure-exploiting ADMM formulation in Algorithm~\ref{alg:ouradmm}, and we show a detailed derivation in  \mbox{Appendix~\ref{app:proofouradmm}}.
We use the Lagrange multipliers~$\lambda_{\zeta}$,~$\lambda_{\epsilon}$, which have the same size and partition as~$y$. Further, we introduce subsystem-individual penalty parameters $\rho_i>0$, a balancing parameter $\beta\in(0,1]$, and the modified projection operation
\begin{IEEEeqnarray}{rCl}
\label{eq:modifiedProjectione}
  {\bar{\proj}_{D\epsilon =d }({\cdot})} = {E_{\epsilon}^{-\nicefrac12}} \proj\nolimits_{D{E_{\epsilon}^{-\nicefrac12}}\epsilon =d }({E_{\epsilon}^{\nicefrac12}} {\cdot}),
  \end{IEEEeqnarray}
where $E_{\epsilon}\sequal (1\sminus\beta ) \operatorname{diag}(\rho_1 I_{\mathrm y_1},\sdots, \rho_M I_{\mathrm y_M})$. In~\circled{2.1}, $\beta$ adjusts the balance between the regularization terms and therefore can affect the convergence speed. We discuss its choice in Section~\ref{sec:parameterChoice}.
 Proposition~\ref{prop:ouradmmconvergence}, proven in \mbox{Appendix~\ref{app:proofouradmm}}, describes how we use Algorithm~\ref{alg:ouradmm}  to solve the original problem~\eqref{eq:mpcprog}.

\begin{algorithm}
\caption{Structure-Exploiting ADMM}\label{alg:ouradmm}
\begin{algorithmic}[1]
\Statex \hspace{-1em}\textbf{repeat}
\algrenewcommand{\alglinenumber}[1]{\hspace{1.4ex}\footnotesize\color{lightblue}\circled{2.1}}
\State  $\forall i\hspace{-.5ex}: \textcolor{lightblue}{y_i}  \leftarrow \arg\hspace{-.3ex}\min\limits_{\textcolor{lightblue}{y_i}} \, \tfrac12 \textcolor{lightblue}{y_i\hspace{-0.55ex}}^\top \hspace{-.5ex} \mathcal Q_i \textcolor{lightblue}{y_i} + q_i^\top \hspace{-.5ex} \textcolor{lightblue}{y_i}
+ 
\frac{{\rho_i}}{2} \Big[ \beta \| \textcolor{lightblue}{y_i} \sminus \textcolor{lightred}{\zeta_i} \sminus \textcolor{lightred}{\lambda_{\zeta_i}} \|_2^2
$ \vspace{-.5ex}

\Statex \hfill$  
+
{(1\sminus\beta)} \| \textcolor{lightblue}{y_i} \sminus \textcolor{lightgreen}{\epsilon_i} \sminus \textcolor{lightgreen}{\lambda_{\epsilon_i} } \|_2^2 \Big]$\vspace{-.5ex}

\Statex \hspace{2.11cm}  $\text{s.t.} \,\, C_i \textcolor{lightblue}{y_i}=c_i$
\vspace{.8ex}
\algrenewcommand{\alglinenumber}[1]{\hspace{1.4ex}\footnotesize \color{lightred}\circled{2.2}}
\State  $ \forall i\hspace{-.5ex}: \textcolor{lightred}{\zeta_i} \leftarrow \proj_{\mathcal Y_i} \left(  \textcolor{lightblue}{y_i}   -  \textcolor{lightred}{\lambda_{\zeta_i}} \right)$ 
          \rlap{\smash{ \raisebox{-1.65ex}{$\bigg\} \hspace{-1mm}\begin{tabular}{l}independent\\execution\end{tabular}$}}}
\vspace{0.8ex}
\algrenewcommand{\alglinenumber}[1]{\hspace{1.4ex}\footnotesize \color{lightgreen}\circled{2.3}}
\State  $ {\epsilon}  \leftarrow { {\bar \proj}_{D{ \textcolor{lightgreen}{\epsilon}} = d} {\left(   \textcolor{lightblue}{y}\, \textcolor{black}{ - } \, \textcolor{lightgreen}{\lambda_\epsilon}  \right)} }$
\vspace{1.2ex}
\algrenewcommand{\alglinenumber}[1]{\hspace{1.4ex}\footnotesize \color{black}\circled{2.4}}
\State  $ \bmat{ \textcolor{lightred}{\lambda_{\zeta}} \\ \textcolor{lightgreen}{\lambda_{\epsilon} } } \leftarrow   \bmat{\textcolor{lightred}{\lambda_{\zeta}} \\ \textcolor{lightgreen}{\lambda_{\epsilon} } }-  \bmat{ \textcolor{lightblue}{y} -\textcolor{lightred}{\zeta} \\ \textcolor{lightblue}{y} - \textcolor{lightgreen}{\epsilon} }$
\end{algorithmic}
\end{algorithm}

\begin{samepage}
\begin{proposition}\hspace{1mm}
\begin{itemize}
	  \item[$(i)$] If $\beta \shortin (0,1)$ and  $\rho_i\hspace{-0.2ex}>\hspace{-0.2ex}0$ for all $i$, Algorithm~\ref{alg:ouradmm} converges to a fixed point $y^\star\hspace{-0.4ex} \sequal \zeta^\star\hspace{-0.4ex}\sequal\epsilon^\star$ that is optimal for~\eqref{eq:mpcprogsplitting} and~\eqref{eq:mpcprog}.
	\item[$(ii)$] If $M\sequal\beta \sequal 1$ and $\rho_1\sequal\rho$, Algorithm~\ref{alg:ouradmm} reduces to Algorithm~\ref{alg:conventionalADMM}. 
\end{itemize}
  \label{prop:ouradmmconvergence}
\end{proposition}
\end{samepage}

In Appendix~\ref{app:proofouradmm}, we introduce $\beta$ and $\rho_i$ through a metric selection technique~\cite{stellato2017osqp,giselsson2017linear}, which makes  Algorithm~\ref{alg:ouradmm} part of the standard ADMM family. The convergence statement in Proposition~\ref{prop:ouradmmconvergence} leaves out the case of $\beta\sequal1$, $M>1$, for which the algorithm is unsuited as it solves~\eqref{eq:mpcprogfinal} without~\eqref{eq:mpcfinal-coupling}.

\begin{table}[t]
\setlength\tabcolsep{6pt} 
\renewcommand{\arraystretch}{1.0}
  \centering
  \caption{Complexities for each step in Algorithm~\ref{alg:ouradmm}} \label{tab:StepComplexities}
  \begin{tabular}{clllll}
    \toprule
      \phantom{\circled{2.1}} & use case & threads & complexity of  the longest thread $\mathcal O(\cdot)$\\
    \midrule
    \circled{2.1} &  & $M$ & $N\max_i\mathrm x_{i}^2$ \\
    \circled{2.2} &  & $2\MN$ & $\max\nolimits_i(\max\{\cost(\proj\nolimits_{\mathcal X_i}), \cost(\proj\nolimits_{\mathcal U_i})\})$ \\
    \circled{2.2} & box & $\MN$ & $\max_i \mathrm x_i$ \\
    \circled{2.3} &  & $N$ & $\mathrm w^2$ \\
    \circled{2.3} & out-1 & $\MN$ &$\max_i \mathrm w_i^2$ \\
    \circled{2.4} &  & $2\MN$ &$\max_{i} \mathrm x_i$ \\
  \end{tabular}
\end{table}
\begin{figure}[t]
  \centering
  \includegraphics[width=1\columnwidth]{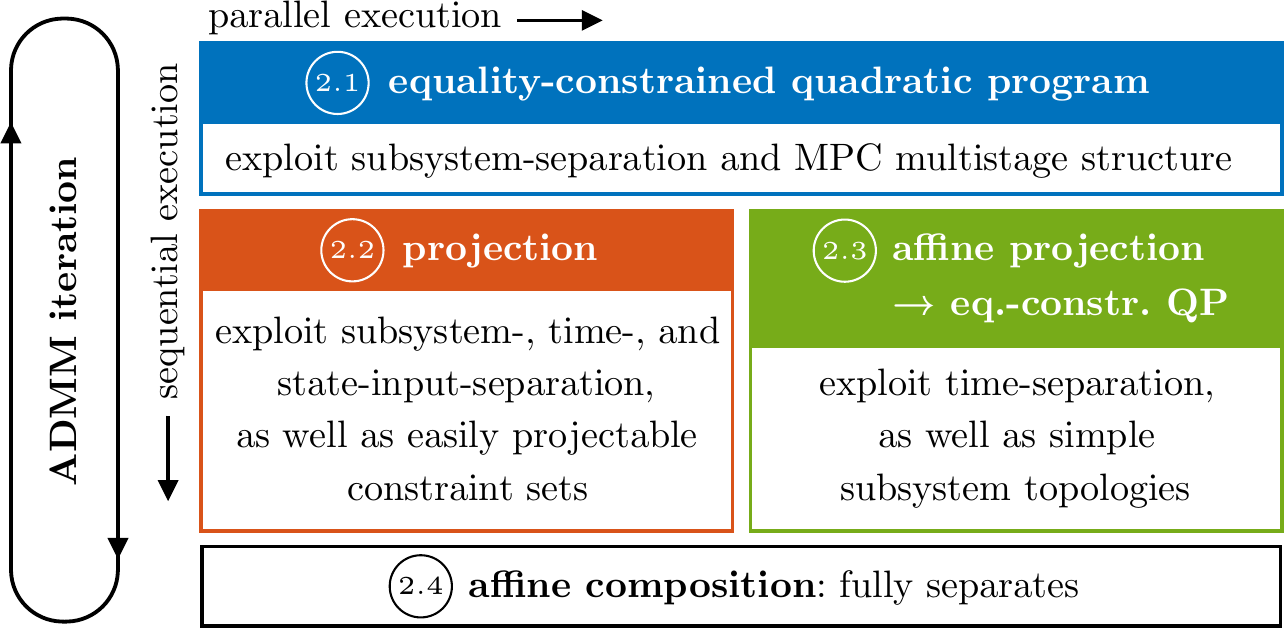}
  \caption{Illustration of the composition of~Algorithm~\ref{alg:ouradmm}. Each  step exploits different types of structure, which improves the computational efficiency.}
  \label{fig:admScheme}
\end{figure}
\begin{table}[!t]
\setlength\tabcolsep{2.5pt} 
\renewcommand{\arraystretch}{1.3}
  \centering
  \caption{Complexity-per-iteration for Algorithm~\ref{alg:ouradmm}} \label{tab:complexities}
  \begin{tabular}{clllll}
    \toprule
       & use case & threads & complexity of the longest thread $\mathcal O(\cdot)$\\
       \midrule
    $\#1$  & $M\sequal1$, box & $1$ &  $N\mathrm x^2 + N \mathrm x$\\
    $\#2$  &  box & $1$ & $N (M\max\nolimits_i \mathrm x_i^2 +  M\max\nolimits_i \mathrm x_i  + \mathrm w^2)$\\
    $\#3$  &  box &  $2\MN$ & $N \max\nolimits_i \mathrm x_i^2 + \max\{ \max\nolimits_i \mathrm x_i, \mathrm w^2\} $\\
    $\#4$ &   box, out-1 &  $1$ & $\MN (\max\nolimits_i \mathrm x_i^2 +  \max\nolimits_i \mathrm x_i + \max\nolimits_i \mathrm w_i^2 )$\\
    $\#5$  &  box, out-1 &  $2\MN$ & $N \max\nolimits_i \mathrm x_i^2 + \max\nolimits_i \{ \mathrm x_i, \mathrm w_i^2 \}$\\
  \end{tabular}
  \renewcommand{\arraystretch}{1}
\end{table}

\subsection{Efficient Implementation and Computational Complexity}
\label{sec:efficientImplementation}
Step~\circled{2.1} is an equality-constrained quadratic program. With the Schur complement method in~\cite[Sec.~16.2]{Nocedal2000}, we obtain the closed-form solution 
	  \begin{IEEEeqnarray}{rCl}
	  \label{eq:step1solvepenalty}
	  y_i &\hspace{-.1ex}\leftarrow\hspace{-.2ex}& {{\mathcal M}_i}\hspace{-.3ex} \big(\sminus \tfrac1{{\rho_i}} q_i  \hspace{-.1ex}+\hspace{-.1ex} {\beta}( \lambda_{\zeta_i} \splus  \zeta_i) \hspace{-.1ex}+\hspace{-.1ex} {(1\sminus\beta)}( \lambda_{\epsilon_i} \splus \epsilon_i)\big) \hspace{-.1ex}+\hspace{-.1ex} {{\mathcal N}_i} c_i, \quad\quad
	  \end{IEEEeqnarray}
	  where $ {\mathcal N}_i =  {\mathcal P}_iC^\top_{i} (C_{i} {\mathcal P}_iC_{i}^\top)^{-1}$,
	  $ {\mathcal M}_i = (I - \mathcal N_i C_{i}) {\mathcal P}_i$, and 
	  ${\mathcal P}_i = {{\rho_i}}(\mathcal Q_i + {{\rho_i}} I)^{-1}$. An efficient implementation of~\eqref{eq:step1solvepenalty} exploits the MPC multistage structure~\cite{domahidi2012efficient}, i.e., it exploits the fact that $C_i$ is banded (see Appendix~\ref{app:StepComplexities} for details).

Similarly, the modified affine projection in~\circled{2.3} can be written as an equality-constrained quadratic program and solved by
	  $
	  \epsilon \hspace{-0.1ex}\leftarrow\hspace{-0.1ex} { {\mathcal D}} (  { y} \hspace{-0.1ex}-\hspace{-0.1ex} { {\lambda}}_\epsilon )\hspace{-0.1ex}+\hspace{-0.1ex}  { {\mathcal E}} d,
	$
	where ${\mathcal E} \sequal {E_{\epsilon}^{-1}} D^\top (D {E_\epsilon^{-1}} D^\top)^{-1} $ and ${\mathcal D} \sequal I - {\mathcal E} D$. Further, we consider a permutation matrix~$P$ that sorts~$y$ for time, i.e., $Py = \bar y = [ \bar y^1;\hdots; \bar y^{N+1}]$ with
$\bar y^1 \sequal [u^1;w^1]$, $\bar y^{k=2,\dots,N} \sequal [x^k;u^k;w^k]$, and $\bar y^{N+1} \sequal x^{N+1}$. In the permuted coordinates, we use $\bar{\mathcal{D}} \sequal  P\mathcal D  P^\top$ and $\bar{\mathcal{E}} \sequal  P\mathcal E$, which are then block-diagonal. Hence, we can solve~\circled{2.3} with
\begin{IEEEeqnarray}{rCl}
\label{eq:step3sol}
  \forall k = 1,\sdots, N\splus1\hspace{-0.5ex}:\quad \,\, {\bar \epsilon}^k \leftarrow \bar{\mathcal D}^k ( \bar { y}^k - {\bar {\lambda}}_\epsilon^k )+  {\bar{\mathcal E}^k} d^k,\quad
\end{IEEEeqnarray}
 where ${\bar {\lambda}}_\epsilon \sequal P\lambda_\epsilon$, $\bar \epsilon \sequal P\epsilon$, and the partition in $k$ matches $\bar y^1, \sdots, \bar y^{N+1}$. A special case is when each subsystem influences at most one other subsystem through a virtual input. In this case, $\bar{\mathcal D}^k$  and $\bar{\mathcal E}^k$ decompose further along $\mathrm  w_1,\sdots,\mathrm w_M$, as the virtual inputs are independent of each other. Consequently,~\eqref{eq:step3sol} can be executed with $\MN$ parallel threads.

\begin{figure}[!t]
  \centering
  \includegraphics[width=.98\columnwidth]{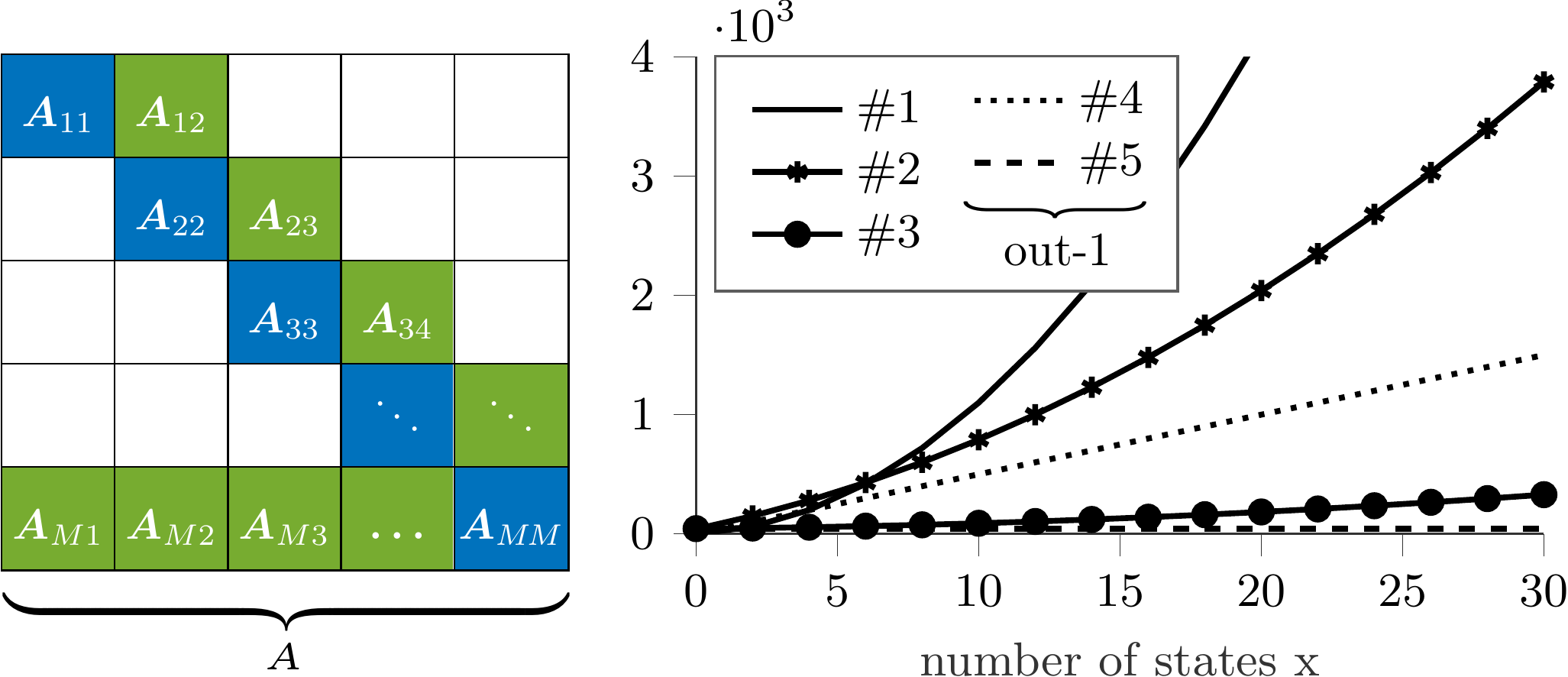}
  \caption{\textit{Left:} Occupation pattern of an example dynamics matrix. \textit{Right:} Complexities-per-iteration in Table~\ref{tab:complexities} for a growing system size.}
  \label{fig:complexityGrowth}
\end{figure}

Table~\ref{tab:StepComplexities} shows the computational complexities of each algorithm step, justifications are shown in Appendix~\ref{app:StepComplexities}. The Landau symbol~$\mathcal O$ describes the order of required scalar multiplications and additions. We denote the computational costs for projecting onto $\mathcal X_i, \mathcal U_i$ with $\cost(\Pi_{\mathcal X_i}), \cost(\Pi_{\mathcal U_i})$. We consider two special cases:  `box' denotes the case where $\mathcal X_i, \mathcal U_i$ are easily-projectable box constraints; and `out-1' denotes the case where each subsystem affects at most one virtual input.
Aside from the structure in the algorithm steps, the composition of Algorithm~\ref{alg:ouradmm} also offers a large potential for parallelization:\looseness=-1
\begin{itemize}
	\item \circled{2.1} and~\circled{2.2} are separate for each subsystem
	\item \circled{2.2} consists of $\MN$ separate projections onto $\mathcal X_i, \mathcal U_i$
	\item \circled{2.3} separates for each step along the prediction interval
	\item \circled{2.2} and \circled{2.3} are independent from each other.
\end{itemize}
	We summarize the composition in Fig.~\ref{fig:admScheme}. 
Table~\ref{tab:complexities} shows the complexity-per-iteration of Algorithm~\ref{alg:ouradmm}. The case $M\sequal1$, independent~of~$\beta$, has the same complexity as conventional ADMM. 
We consider an example system to illustrate the complexities in Fig.~\ref{fig:complexityGrowth}. In Table~\ref{tab:complexities}, we see that the obtained values follow from the problem parameters $M$, $N$, $\{\mathrm x_i, \mathrm w_i\}$, the number of threads, and the use case.
We use $A_{ij} \in \mathbb R^{2\times2}$ for all~$i,j$. We assume that  each $A_{i,i+1}$ has rank~$1$, i.e.,  $\mathrm w_i \sequal1$ for $i\neq M$. For the last component, we assume $\mathrm w_M \sequal  2$. 
We use $N\sequal 10$ and we consider systems with~\mbox{$1$-$15$} diagonal components, i.e., $\mathrm x \in  \{2,4,6,\sdots,30 \}$. In Fig.~\ref{fig:complexityGrowth} on the right, we show the corresponding complexities.
 The cases~$\#4$,$5$ only apply if $A_{i,i+1}\sequal 0$ or $A_{M,i}\sequal 0$, i.e., the case `out-1' applies.
We see that structure exploitation (\#2-5) compares favorably to conventional ADMM (\#1), even for single-thread implementations (\#2,4).
We also see that the benefit from structure exploitation grows with the system size. As presented in Fig.~\ref{fig:setup}, the complexity-per-iteration only defines a part of the overall execution cost. We analyze the remaining part, namely the required number of iterations, in Section~\ref{sec:separationTendency}.

\subsection{Parameter Choice}
\label{sec:parameterChoice}

We denote the smallest and largest eigenvalues of a matrix with $\operatorname{eig}_{min}, \operatorname{eig}_{max}$. We consider an orthonormal null space basis for \mbox{$C_i\in\mathbb R^{N\mathrm x_i \times \mathrm y_i}$}, and use the basis vectors as columns in  \mbox{$Z_i\in\mathbb R^{\mathrm y_i\times N\mathrm x_i}$}, which leads to $C_iZ_i\sequal0$ and $Z_i^\top Z_i \sequal I_{N\mathrm x_i}$. 

\begin{proposition}
 We assume that $Z_i^\top \mathcal Q_i Z_i$ is positive definite and $\beta \in (0,1)$. 
  The optimal penalty parameters for improving the worst-case convergence rate of Algorithm~\ref{alg:ouradmm} are
  \begin{IEEEeqnarray}{rCl}
  \label{eq:optimalrho}
    \rho_i^\star = \sqrt{\operatorname{eig}_{min}(Z_i^\top \mathcal Q_i Z_i)\operatorname{eig}_{max}(Z_i^\top \mathcal Q_i Z_i)},
  \end{IEEEeqnarray}
  \label{prop:optimal param}
  where $i=1,\sdots,M$.
\end{proposition}

\ifArxiv
We show a proof of Proposition~\ref{prop:optimal param} in Appendix~\ref{app:proofoptimalparam}, and we provide additional details in the ancillary material. 
\else
We show a proof of Proposition~\ref{prop:optimal param} in Appendix~\ref{app:proofoptimalparam}, and we provide additional details and technical considerations in~\cite{Rey2018}. 
\fi
The proof also shows that $\beta$ is canceled out when we derive~$\rho_i^\star$, i.e., the parameters can be chosen independently. Proposition~\ref{prop:optimal param} suggests that individual penalty parameters are indeed useful, as their optimal choice is different from making them all the same. 
The optimal parameters~\eqref{eq:optimalrho} are valid for any quadratic program of type~\eqref{eq:mpcprogfinal}. For MPC problems in particular, 
and similar as for comparable results~\cite{raghunathan2014admm,raghunathan2014optimal,giselsson2017linear}, we observe that we often can improve the practical performance further by increasing the penalties above $\rho_i^\star$, which places an additional weight on the regularization terms in \circled{2.1}.

As noted before, $\beta\in(0,1]$ adjusts the regularization balance in~\circled{2.1}. The two regularization terms are equally weighted for~$\beta \sequal \tfrac12$. By increasing~$\beta$ we emphasize the influence of $(\zeta_i,\lambda_{\zeta_i})$ over $(\epsilon_i,\lambda_{\epsilon_i})$, and vice versa. We can also adapt~$\beta$ during the algorithm iteration, similar to $\rho$ in~\cite[Sec.~3.4.1]{boyd2011distributed}. 
However, our numerical results suggest that the effect of $\beta$ is small for the cases considered below. For this reason, we restrict our attention to the nominal values $\beta \in \{\tfrac12,1\}$ in the sequel. For $\beta \sequal \tfrac12$ and $M>1$, we speak of structure-exploiting ADMM. With $\beta\sequal M\sequal1$, we refer to conventional ADMM.
 
\section{Separation Tendency}
\label{sec:separationTendency}
We derive a quantitative measure of system structure, called the separation tendency. 
We use the separation tendency as a heuristic indicator for the required number of algorithm iterations relative to conventional ADMM. 
In combination with the complexity results in Table~\ref{tab:complexities}, this makes it possible to anticipate the execution cost of structure exploiting ADMM. 

\subsection{System Flow and Link Usage} Our goal is to quantify the interaction between system components. We first consider the unpartitioned $\Delta$-system
\begin{IEEEeqnarray}{rCL}
  \Delta x^{k+1} = A \Delta x^k + B \Delta u^k,
  \label{eq:detlasystem}
\end{IEEEeqnarray}
where $\Delta x^k \sequal x^k - x^{k-1}$, $\Delta u^k \sequal u^k - u^{k-1}$, $k\geq0$, and we use the convention $(x^{-1}, u^{-1}) \sequal (0,0)$. The $\Delta$-system describes the changes in the original system~\eqref{eq:initsys}. 

\begin{definition}
\label{def:flow}
  The \textit{system flow} $\Phi^k = [\Phi^k_A, \Phi^k_B] \in \mathbb R^{{\mathrm x}\times ({\mathrm x + \mathrm u})}$ is composed of the state-to-state and input-to-state flow
  \begin{IEEEeqnarray}{rCll}
  \IEEEyesnumber\eqlabel{eq:flowdef}\IEEEyessubnumber*
    \Phi^k_A &=& A\operatorname{diag}\left(\Delta x^k \right) &\in \mathbb R^{\mathrm x \times x}\\
    \Phi^k_B &=& B\operatorname{diag}\left(\Delta u^k \right)\,&\in \mathbb R^{\mathrm x \times u}.
   \end{IEEEeqnarray}
\end{definition}

We use the operation  \mbox{$A\diag (\Delta x)\in \mathbb R^{\mathrm x\times \mathrm x}$} to analyze system-internal effects.  With \mbox{$A\diag (\Delta x) 1_{\mathrm x \times 1} \sequal A\Delta x$}, it becomes clear that we can understand $A\diag (\Delta x)$ as an intermediate step that leads to the matrix-vector product $A\Delta x$. By following this relation, we obtain
\begin{IEEEeqnarray}{rCl}
\label{eq:summedFlowDrivingStates}
  x^{k+1} = x^k +   \Phi^{k} 1_{(\mathrm x + \mathrm u) \times 1},
   \end{IEEEeqnarray}
which clarifies that the system flow $\Phi^{k}$ describes the state transition. A central characteristic of $\Phi^{k}$ is that it details the transition contribution for each of the $\mathrm x^2$ state-to-state links and $\mathrm x \cdot\mathrm u$ input-to-state links. 
The link usage, which we introduce next,  measures the flow that moves through each system link in response to a  unit input impulse~$\delta^{k}$ with $\delta^{k=0}\sequal1$ and $\delta^{k\neq 0}\sequal0$.

\begin{definition}
\begin{samepage}
\label{def:linkusage}
  The \textit{link usage} $\Gamma = [\Gamma_A, \Gamma_B]\in \mathbb R^{{\mathrm x}\times ({\mathrm x + \mathrm u})}$ is assembled element-wise with
  \begin{IEEEeqnarray}{rCl}
  \label{eq:linkusagedef}
  \Gamma_{ij} &=& \left(\sum\nolimits_{k=0}^\infty {|\Phi_{ij}^k|}^2\right)^{\nicefrac12} \in \mathbb R,
   \end{IEEEeqnarray}
   where $\Phi^k \sequal \left[\Phi_{ij}^k\right]$, $\Gamma \sequal \left[\Gamma_{ij}\right]$, and the sequence~$\{\Phi^k\}$ results from \mbox{$u^k = \delta^k$} and $x^0=0$.
   \end{samepage}
\end{definition}

The link usage $\Gamma$ analyzes the system flow over time by using an element-wise $\mathcal L_2$ norm~\cite[Chapter~2]{okuyama2014discrete}. If $\Gamma_{ij}$ is large, then the respective system link is used intensively.

\subsection{Separation Tendency}
The separation tendency compares the link usage for internal and external elements, which provides a relative measure for the concentration of flow inside and outside of virtual subsystems. In contrast to the previous concepts, the separation tendency depends on the system $\textit{and}$ its partition.

\begin{definition}\label{def:s}
  The \textit{separation tendency} $s\in\mathbb R$ is defined by
  \begin{IEEEeqnarray}{rCl} 
  \IEEEyesnumber\eqlabel{eq:defs} \IEEEyessubnumber*
    s_i &=& \frac{\tfrac{1}{\operatorname{\#int}_i}\sum_j \internal{\Gamma}_{ij}}{\tfrac{1}{\operatorname{\#int}_i}\sum_j \internal{\Gamma}_{ij} + \tfrac{1}{\operatorname{\#ext}_i}\sum_j  \external{\Gamma}_{ij}} \label{eq:defsa}
   \\
   s &=&  \tfrac{1}{\mathrm x} \sum_{i=1}^{\mathrm x}s_i
   ,
   \end{IEEEeqnarray}
   where $\Gamma \sequal \internal{\Gamma} + \external{\Gamma} \sequal [\internal{\Gamma}_A, \internal{\Gamma}_B] + [\external{\Gamma}_A, \external{\Gamma}_B]$ is an internal-external decomposition, and $\operatorname{\#int}_i$, $\operatorname{\#ext}_i$ are the numbers of internal and external elements in the $i$-th row of~$\Gamma$.
   \end{definition}

 If the separation tendency $s$ is large, then internal links predominantly influence the states, which signals a clear subsystem separation. If the separation tendency is small, then the system states are dominated by external flow, which signals that the chosen partition is unsuited. 

\begin{samepage}
\begin{proposition}\hspace{1mm}
\label{prop:stproperties}
  \begin{itemize}
  	\item[$(i)$] $s$ exists if the system $(A,B)$ is controllable, and $A$ is semi-convergent, i.e., $\lim\nolimits_{k\rightarrow \infty} A^k$ exists.
    \item[$(ii)$] If $s$ exists, then $0\leq s \leq 1$.
    \item[$(iii)$] $s$ is invariant to diagonal state and input transformations.
  \end{itemize}
\end{proposition}
\end{samepage}

We show a proof of Proposition~\ref{prop:stproperties} in Appendix~\ref{app:proofstproperties}. We use semi-convergence, which is a weaker condition than asymptotic stability, but stronger than marginal stability. The existence of the separation tendency~$s$ implies that the link usage~$\Gamma$ is finite and does not contain zero rows, which prevents that the denominator in~\eqref{eq:defsa} becomes zero.  
Property $(iii)$ makes clear that~$s$ is unaffected by state and input transformations $[x;u]=T[\bar x;\bar u]$, where $T$ is an invertible and diagonal matrix. This property is important as it makes the separation tendency invariant to simple state and input scaling. 

\begin{figure}[t]
  \centering
  \includegraphics[width=.95\columnwidth]{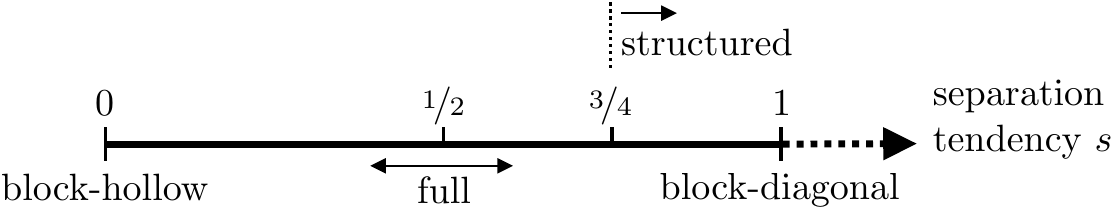}
  \caption{Empirical characterization of the separation tendency $s\in [0,1]$.}
  \label{fig:tendencyRange}
\end{figure}

 It is easy to show that $s\sequal 0$ for block-hollow systems and $s\sequal 1$ for block-diagonal systems. Full systems are placed in-between with $s\approx 0.5$. Fig.~\ref{fig:tendencyRange} illustrates the range of $s$.  We call a system structured if it can be partitioned with a large separation tendency. An empirically reasonable threshold between structured and unstructured systems is $s\sequal \tfrac34$.


\begin{example}
\label{ex:example2}
We consider the full system
\begin{IEEEeqnarray}{lCl}
 \label{eq:examplesystem2}
 x^{k+1} = \underbrace{\bmat{\nicefrac12 & \nicefrac12 \\ \nicefrac12 & \nicefrac12}}_{A} x^{k} + \underbrace{\bmat{1 \\ 1}}_B u^{k}
  \end{IEEEeqnarray}
and a partition with $\{\mathrm x_i\} \sequal \{1,1 \}$ and $\{\mathrm u_i\} \sequal \{1,0 \}$. We observe that the system matrix is semi-convergent (it even is idempotent, i.e., $(A)^k=A$). The system is not controllable, therefore $s$ is not guaranteed to exist a-priori.  For an input impulse~$\delta^{k}$, the system flow and the link usage become
\newcommand\cellb{\cellcolor{lightblue!100}\color{white}}
\begin{IEEEeqnarray*}{rCl}
  \Phi^{k} &=&\delta^{k}\tfrac12\bmat{0 & 0 & 1\\0 & 0 & 1} + \delta^{k-1}\tfrac12\bmat{1 & 1 & -1\\1 & 1 & -1}\\
  \Gamma &=& \Bigg[\begin{array}{cccc} \cellb \nicefrac12 & \nicefrac12 & \cellb \sqrt{2}\\\nicefrac12 & \cellb \nicefrac12 & \sqrt{2}\end{array}\Bigg],
 \end{IEEEeqnarray*} 
 where we shade the link usage according to the partition. As all elements in $\Gamma$ are non-zero and finite, the separation tendency exists. We obtain $s=\tfrac12$,
 which suggests that the considered partitioned system is not suited for structure exploitation.
\end{example}

\subsection{Algorithm Performance Indication}
\label{sec:empiricalResults}
\begin{figure}[t]
\setlength\fwidth{0.871\columnwidth}
  \centering
  \input{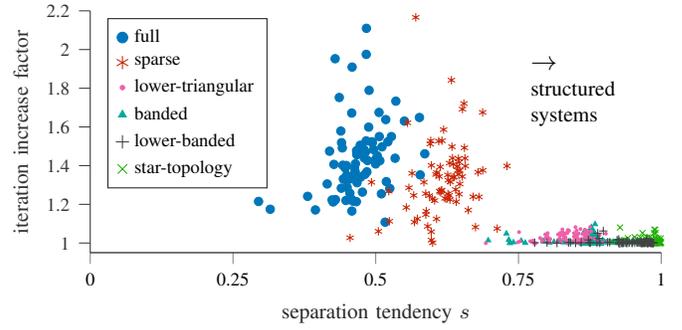}
\caption{Iteration increase over the separation tendency. Each of the $480$ dots represents an example system. For each system, we generate twenty problem instances, run both algorithms, and depict the average increase.}
\label{fig:iterationComparison}
\end{figure}

We show that $s$ is an indicator for the performance of system structure exploitation. More specifically, we collect empirical evidence that~$s$ relates to the growth in required algorithm iterations when we switch from conventional to structure-exploiting ADMM. 
We consider systems in six categories: full, sparse, lower-triangular, banded, lower-banded, and star-topology; 
\ifArxiv
as described in the ancillary material.
\else
as described in~\cite{Rey2018}.
\fi
Additionally, we consider the dimensions \mbox{$\mathrm x \in \{5, 10, 20, 40\}$}, which result in~$24$ combinations. For each combination, we generate twenty pairs of system matrices, which leads to a test set of~$480$ systems.
   For each system, we then generate twenty problems of type~\eqref{eq:mpcprog}. For simplicity, and as $s$ only depends on the partitioned system, we set $\mathcal X\sequal \mathbb R^{\mathrm x}$ $\mathcal U\sequal \mathbb R^{\mathrm u}$. We solve the final $9600$ problems with conventional and structure-exploiting ADMM. For the structure-exploiting case, we choose a partition that fits to the problem type and dimension 
\ifArxiv
as described in the ancillary material.
\else
as described in~\cite{Rey2018}.
\fi
   For the penalty parameters, we set $\rho\sequal\rho_i\sequal1$. 
We measure the number of iterations that are necessary to converge within a certain accuracy of a precomputed solution. We then compute the iteration increase factor when we switch between algorithms, and we average this factor over each system's twenty initial conditions. Fig.~\ref{fig:iterationComparison} illustrates the result. The key observation is that a large separation tendency indicates a low iteration increase. Combined with the complexity results in Table~\ref{tab:complexities}, this allows us to assess the algorithm performance as defined in Fig.~\ref{fig:setup}, particularly without having to implement and benchmark the algorithm first. 

\subsection{Related Concepts and System Properties}
Several quantities that are used in the literature are related to~$\Phi$,~$\Gamma$, and~$s$. Spectral clustering methods~\cite{ageev2000approximation,frieze1997improved,hespanha2004efficient} interpret the system as a weighted graph, where the elements of~$[A,B]$ determine the edge weights. These methods can be used to determine a system partition by cutting possibly low-weighted edges. In contrast to the static weights~$[A,B]$, we use impulse-response-based dynamical links for~$\Gamma$. 
We can assess the value of this dynamical concept by redrawing Fig.~\ref{fig:iterationComparison} while we use~$[A,B]$ instead of $\Gamma$ to obtain~$s$. We then observe that this static version of~$s$ is significantly less indicative for the algorithm iteration growth.
Furthermore, clustering methods are sensitive to diagonal state and input transformations.

The similarity of~$A,B$ to block-diagonal matrices can be an intuitive structure measure as it directly relates to the sparsity of $\external{A},\external{B}$. 
In~\cite{alyani2017diagonality}, a range of diagonality measures is discussed; generalizations to block-diagonality are straightforward. 
Same as before, these measures ignore dynamic interaction, are less indicative for the iteration growth, and are sensitive to diagonal state and input transformations.

The computation of $\Gamma$ is similar to the computation of 
 the $\mathcal H_2$ system norm, where the $\mathcal L_2$ norm is applied to the system's impulse response~\cite[Eq.~(2.37)]{toscano2013structured}, \cite[Eq.~(2.167)]{oppenheim1999discrete}. 
 Two main differences separate the concepts. First, $\Gamma$ is based on the \mbox{$\Delta$-system}, which makes it finite for a wider range of cases. Second, $\Gamma$ is matrix-valued, which underlines the focus on the system-internal state-to-state and input-to-state links.
In contrast, the $\mathcal H_2$ norm is scalar-valued, even if we use the full state vector as system output~\cite[Eq.~(2.37)]{toscano2013structured}.

Another related concept is the balanced realization \cite[Sec.~4.2]{gu2012discrete}, which is a state space system representation with identical and diagonal controllability and observability Gramians. The diagonal elements then quantify the influence of each state on the input-output behavior. 
The concept resembles~$s$ in its quantitative description of system-internal relations. However,~$s$ is based on $\Phi$ with $\mathrm x^2 + \mathrm x\mathrm u$ elements, while the Gramians only have $\mathrm x$ diagonal elements. Hence, $\Phi$ analyzes the system with a higher resolution. Aside from that, another distinction is that $\Phi$ (and therefore $s$) focuses on the state transition, while the Gramians relate to the input-output behavior.

\section{Simulation Study}
We apply  structure-exploiting ADMM to a cascade system, where our method leads to a substantial benefit. We also show a negative example of an unstructured system. More positive examples can be found in~{\cite{Rey2017b,Rey2017a}}.

\subsection{Cascade System}
\label{sec:simulation:cascade}

A cascade system~\cite{cantoni2017structured} is characterized by a lower block-banded dynamics matrix and a block-diagonal input matrix
\begin{IEEEeqnarray}{rCl}
\arraycolsep=.5pt
 x^{k+1} \sequal \underbrace{\bmat{  A_{11} \\ A_{21} & A_{22} \\ & \ddots & \ddots \\ && A_{S(S-1)} & A_{\twoS} }}_{A} x^k \splus
 \underbrace{\bmat{B_{11}\\ & B_{22} \\ & & \ddots \\ & & & B_{\twoS}}}_{B}u^k, \quad\nonumber\\[-7ex]\label{eq:cascadeSystem} \\\nonumber
\end{IEEEeqnarray}
where we use the stages~$i\sequal1,\sdots,S$, each with a dynamics matrix~$A_{ii}$, input matrix~$B_{ii}$, and coupling matrix~$A_{i(i-1)}$. 
Cascade systems are used for irrigation and drainage networks~\cite{li2005water,soltanian2015decentralized}, hydro-power systems~\cite{labadie2004optimal}, and vehicle platoon control~\cite{guo2011hierarchical}. In~\cite{cantoni2017structured}, an interior point method is developed, where the iteration complexity scales linearly in~$S$ and cubically in~$N$. 
For structure-exploiting ADMM, we associate each stage to a virtual subsystem, i.e., $M\sequal S$. Due to the simple subsystem topology, the complexity results~$\#4,5$ in Table~\ref{tab:complexities} apply. Hence, for a single-thread implementation, the iteration complexity scales linearly in~$S$ and~$N$. Further, when we use parallel computation, the complexity becomes constant in~$S$.

 We consider $S\sequal 20$ stages, each with $\mathrm x_i\sequal6$ states and \mbox{$\mathrm u_i\sequal1$} input, resulting in a cascade system with $\mathrm x \sequal 120$ states and $\mathrm u \sequal 20$ inputs.
The system matrices are randomly generated as 
\ifArxiv
described in the ancillary material.
\else
in~\cite{Rey2018}. 
\fi
The stage coupling~$A_{i(i-1)}$ has a rank equal to~$1$. We use~\eqref{eq:mpcprog} with $N\sequal 5$ and box constraints. 
Table~\ref{tab:costresult}  shows estimates for the cost-per-iteration of Algorithm~\ref{alg:ouradmm} in different situations.
\begin{table}[!t]
\setlength\tabcolsep{6pt} 
  \centering
  \caption{Cost-per-iteration of Algorithm~\ref{alg:ouradmm} for the cascade system~\eqref{eq:cascadeSystem}} \label{tab:costresult}
  \begin{tabular}{r@{\hspace{4pt}}lccccc}
    \toprule
      &ADMM-type & $M$  & threads & cost &$\nicefrac{\text{cost}}{(i)}$\\
    \midrule
    $(i)$ & conventional  & $1$  & $1$ & $277676$  &$100\%$\\
    $(ii)$ & structure-exploiting &  $S$  & $1$ & $50220$ &$18.09\%$\\
    $(iii)$ & structure-exploiting & $S$  & $2\MN$ & $3125$  &$1.13\%$\\
  \end{tabular}
\end{table}
In contrast to the analytical bounds in Table~\ref{tab:complexities}, we obtain the computational costs by counting the scalar additions and multiplications in an actual implementation.
This counting strategy is more precise than complexity bounds and takes the remaining matrix sparsity into account. As opposed to timing measurements, it is also less hardware-dependent. For the parallel implementation $(iii)$, we count the operations in the longest thread. We do not account for memory access or data exchange operations, which is justified for many FPGA-type implementations where such operations can be hard-coded.
We see that structure exploitation significantly reduces the computational cost.\looseness=-1

\begin{figure}[!t]
  \begin{tikzpicture}[
  myrect/.style={
    rectangle,
    draw,
    inner sep=0pt,
    fit=#1}
  ]

  \clip (-4.5,-2.95) rectangle (4.2, 2.45);
    \node (img)  { \includegraphics[trim=0 0 0 0,clip,width=0.94\columnwidth]{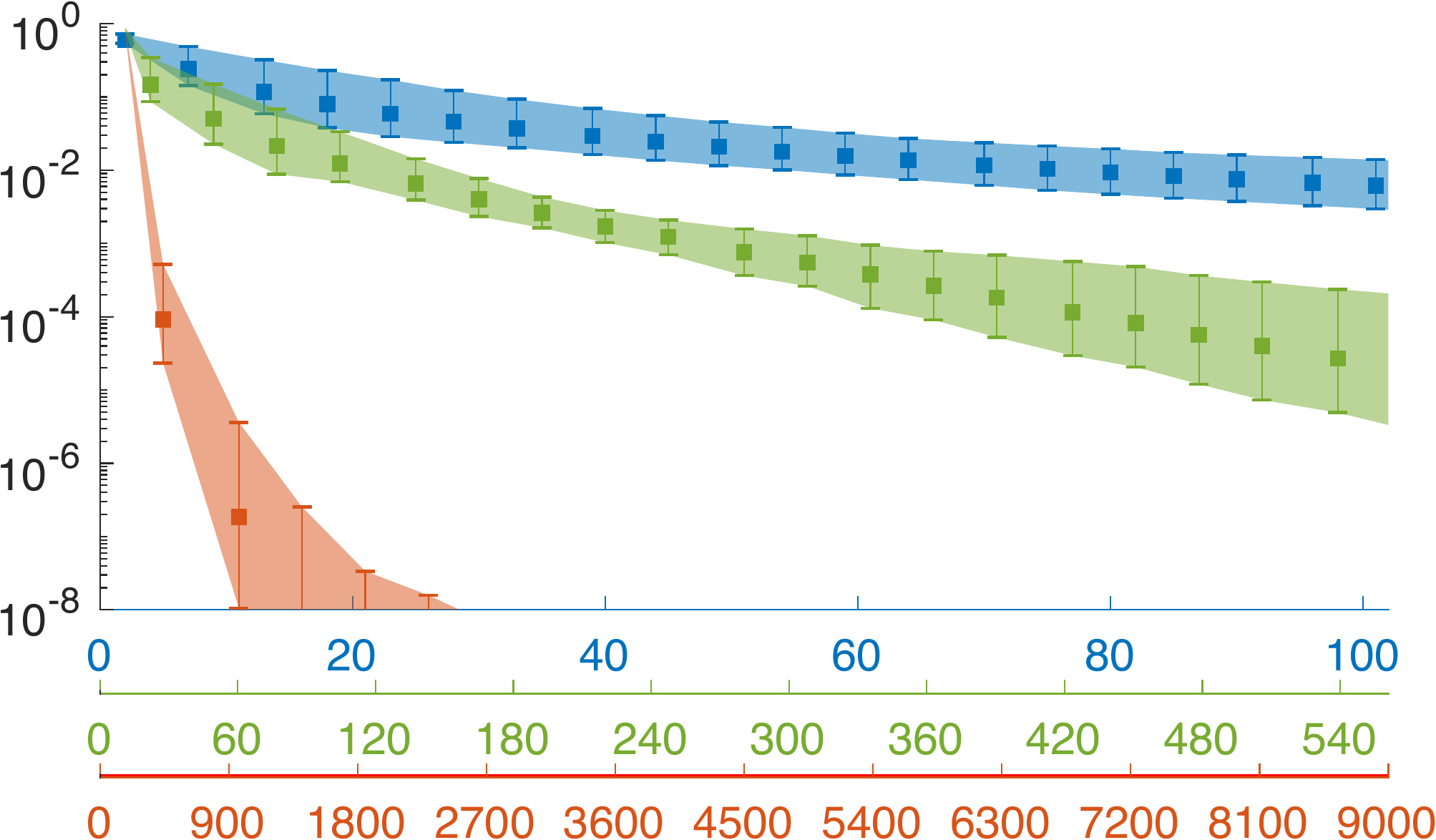}};
    \node[below of = img, node distance=3.75cm, yshift=.95cm] {\footnotesize number of algorithm iterations};
    \node[left of = img, node distance=5.15cm, rotate=90, anchor=center,yshift=-0.8cm] {\footnotesize solution distance $\operatorname{dist}(x,u)$};

    \coordinate (legendboxCoordinate) at (.6,-0.8);
    \node[  name=legendbox,
              draw=white,
              fill=white,
              rectangle,  
              text width=5cm,
              minimum height=0.4cm]   at ($(legendboxCoordinate) + (0.75, 0)$) {};

      \coordinate (legendLine1L) at ($(legendboxCoordinate) + (-2, 0.8)$);
    \coordinate (legendLine2L) at ($(legendboxCoordinate) + (-2, 0.35)$);
    \coordinate (legendLine3L) at ($(legendboxCoordinate) + (-2, 0.05)$);

\coordinate (A) at ($(legendLine1L) + (0, 0.07)$);
\coordinate (B) at ($(legendLine1L) + (0.35, -0.08)$);
\node[myrect={(A) (B)},fill=lightblue, draw=white, line width=0pt] {}; 

\coordinate (A) at ($(legendLine2L) + (0, 0.07)$);
\coordinate (B) at ($(legendLine2L) + (0.35, -0.08)$);
\node[myrect={(A) (B)},fill=lightgreen, draw=white, line width=0pt] {}; 

\coordinate (A) at ($(legendLine3L) + (0, 0.07)$);
\coordinate (B) at ($(legendLine3L) + (0.35, -0.08)$);
\node[myrect={(A) (B)},fill=lightred, draw=white, line width=0pt] {}; 

    \node[anchor=west] at ($(legendLine1L) + (0.35, 0)$) {\footnotesize $(i)$ conventional ADMM};
    \node[anchor=west] at ($(legendLine2L) + (0.35, 0)$) {\footnotesize $(ii)$ \,\,single-thread};
    \node[anchor=west] at ($(legendLine3L) + (0.35, 0)$) {\footnotesize $(iii)$ $2\MN$-thread};

    \node[anchor=west] at ($(legendLine3L) + (2.6, 0.15)$) {$\Big\}$};
      \node[anchor=west] at ($(legendLine2L) + (2.9, 0)$) {\footnotesize structure-exploiting};
      \node[anchor=west] at ($(legendLine3L) + (2.9, 0)$) {\footnotesize ADMM};

     \end{tikzpicture}
   \caption{Convergence of Algorithm~\ref{alg:ouradmm} for the cascade system~\eqref{eq:cascadeSystem}. The markers indicate the geometric mean over $200$ control scenarios; the illustrated range includes $80\%$ of all scenarios, excluding the best and~worst~$10\%$.}
   \label{fig:cascadeConvergence}
 \end{figure}
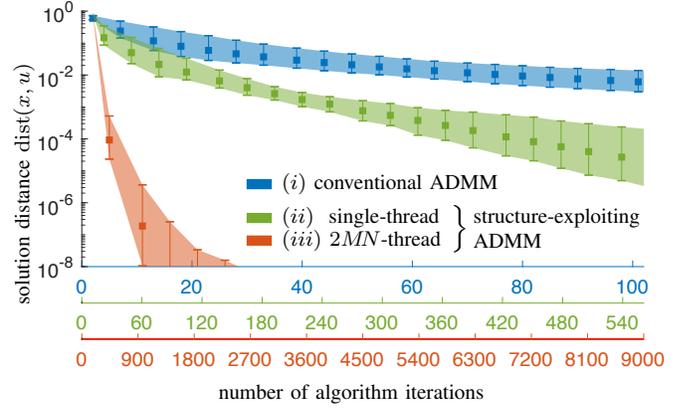

By using the methods in Section~\ref{sec:separationTendency}, we obtain $s\sequal0.975$ for the cascade system.
 We assess the required number of algorithm iterations by applying $(i)$\,-\,$(iii)$  in a range of control situations. More precisely, we generate $200$ feasible instances of~\eqref{eq:mpcprog} 
\ifArxiv
as described in the ancillary material.
\else
 as described in~\cite{Rey2018}. 
\fi
 We compute $\rho_i^\star$ for each subsystem, and we increase the penalty parameters with a factor of~$90$ for improving the MPC performance as discussed in Section~\ref{sec:parameterChoice}.
We then analyze the convergence to a precomputed nonzero solution~\mbox{$(x^\star, u^\star)$} with
\begin{IEEEeqnarray}{rCl}
\operatorname{dist}(x,u)= {
\left\|[x; u] - [x^\star; u^\star] \right\|_2^2
}\,/\,{ 
\left\|[x^\star; u^\star] \right\|_2^2
},
\end{IEEEeqnarray}
where \mbox{$(x, u)$} is the current estimate, extracted after each iteration. We rely on free licenses for Yalmip~\cite{lofberg2004yalmip} and Gurobi~\cite{gurobi}. In  Fig.~\ref{fig:cascadeConvergence}, we show the growing solution accuracy with the number of performed iterations as a statistic over the $200$ problem instances.
 We compare the overall performance of $(i)$\,-\,$(iii)$ by using three different horizontal axes, scaled with the computational costs in Table~\ref{tab:costresult}. For a given point along these axes, we see how many iterations each method can perform while they use the same number of sequential scalar operations.
 We observe that for the cascade system, the exploitation of system structure results in faster convergence, even for a single-thread implementation. 
 When we use the full parallelization potential, high performance is possible. By sequentially combining previously parallel threads, it is also possible to obtain implementations that perform between~$(ii)$ and~$(iii)$ with less than $2\MN$~threads.

\subsection{Unstructured System}
We consider the system from  Example~\ref{ex:example2} with $\{\mathrm x_i\}\sequal\{1,1\}$ and $\{\mathrm u_i\}\sequal\{1,0\}$, i.e.,
\newcommand\cellg{\cellcolor{lightgreen!100}\color{white}}
\newcommand\cellb{\cellcolor{lightblue!100}\color{white}}
\begin{IEEEeqnarray}{rCl}
\label{eq:unstructuredsystem}
  x^{k+1} = \Bigg[\begin{array}{cc} \cellb \nicefrac12 &  \nicefrac12 \\  \nicefrac12 & \cellb \nicefrac12 \end{array}\Bigg] x^k 
  + \Bigg[\begin{array}{cc} \cellb 1 \\  1 \end{array}\Bigg] u^k.
\end{IEEEeqnarray}
The partition is unsuited as it does not align with any visible system structure. This observation is reflected in $s\sequal\tfrac12$. We embed~\eqref{eq:unstructuredsystem} in an MPC setting with diagonal objectives and box constraints, which makes the partition admissible. Table~\ref{tab:costresultUnstructured} shows the numerical cost analysis.
\begin{table}[!t]
\setlength\tabcolsep{6pt} 
  \centering
  \caption[Cost-per-iteration of Algorithm~\ref{alg:ouradmm} for the unstructured system~\eqref{eq:unstructuredsystem}]{\tabular[t]{cc}Cost-per-iteration of Algorithm~\ref{alg:ouradmm} \\ for the unstructured system~\eqref{eq:unstructuredsystem}\endtabular}\label{tab:costresultUnstructured}
  \begin{tabular}{r@{\hspace{4pt}}lccccc}
    \toprule
      &ADMM-type & $M$ & threads & cost &$\nicefrac{\text{cost}}{(i)}$\\
    \midrule
    $(i)$ & conventional  & $1$  & $1$ & $878$  &$100\%$\\
    $(ii)$ & structure-exploiting &  $S$  & $1$ & $1329$ &$151\%$\\
    $(iii)$ & structure-exploiting & $S$  & $2\MN$ & $552$  &$63\%$\\
  \end{tabular}
\end{table}
We see that when we switch from $(i)$ to $(ii)$, the computational cost increases. 
For a parallel implementation $(iii)$, the cost reduces again. However, further simulations show that the overall performance remains  worse than for conventional ADMM. Hence,  structure exploitation only performs well if the controlled system has structure.

\section{Conclusions and Future Work}
We adapt ADMM to exploit structure in MPC. If the controlled system is sufficiently structured, the resulting algorithm scales well, can be specialized with multiple penalty parameters, is highly parallelizable, and shows improved overall performance. Our algorithm reduces the cost-per-iteration, especially for large and structured systems. The cost reduction comes with an increase in necessary algorithm iterations due to the virtual decomposition of the system. We introduce the separation tendency, a  measure of subsystem independence, to relate the iteration increase to the level of structure in the controlled system.
 Finally,  we show a cascade system example, where our structure-exploiting method significantly outperforms conventional ADMM. 
In future work, the separation tendency can be utilized to determine a partition in the first place. Also, a procedure can be developed that finds a state and input transformation for maximizing exploitable structure, while maintaining the partition admissibility. Furthermore, the concepts of system flow and link usage can be used in other areas of system analysis.


%

\appendices
\section{Details and Proofs}
\label{app:proofs}
\subsection{Stacked Problem Formulation}
\label{app:stacking}
For the objective in~\eqref{eq:mpcprogfinal}, we use $\mathcal Q_i \sequal I_{N} \otimes \diag( R_i, 0_{\mathrm w_i},  Q_i)$, $r_{y_i} \sequal \big[r^1_{u_i}; 0_{\mathrm w_i \times 1}; r^2_{x_i} ; \sdots ;r^N_{u_i}; 0_{\mathrm w_i \times 1}; r^{N+1}_{x_i}\big]$, $q_i \sequal - \mathcal Q_i r_{y_i}$, and $K_i = \tfrac12 r_{y_i} \mathcal Q_i r_{y_i}$, where $\otimes$ is the Kronecker product.
We use
\renewcommand*{\arraystretch}{.9}
\begin{IEEEeqnarray*}{l}
\arraycolsep=2.2pt
  \left[C_i, c_i\right] = \hspace{-1mm}  
  \left[ \arraycolsep=1.8pt \begin{array}{ccccccccc|c}  
   B_{ii} &  W_i & -I &  &   & & &&&  - A_{ii} x_{i}^1\\ 
    &   &  A_{ii}&  B_{ii} &  W_i & -I     & & && 0  \\[-1ex]  
  &  & &  & & \ddots     &   & & & \vdots   \\[.2ex]  
  &    &   &   & &  A_{ii} &  B_{ii} &  W_i & -I   & 0\end{array}\right],
  \end{IEEEeqnarray*}
\renewcommand*{\arraystretch}{1}%
  where we understand $C_i\in\mathbb R^{N\mathrm x_i \times \mathrm y_i}$ such that $A_{ii}$ is always below $-I$. For the individual constraints, we choose $\mathcal Y_i$ such that $(u^k_i,w^k_i,x^{k+1}_i) \in \mathcal U_i \hspace{-0.1ex}\times\hspace{-0.1ex} \mathbb R^{\mathrm w_i} \hspace{-0.1ex}\times\hspace{-0.1ex} \mathcal X_i$ for all $k$. For the coupling, we recognize that $[\external{A}, \external{B}, W]$ has a row rank defect of $\mathrm x \hspace{-0.2ex}-\hspace{-0.2ex} \mathrm w$, and we obtain the reduced form \mbox{$[\external{A}_r, \external{B}_r, W_r]\in\mathbb R^{\mathrm w \times (\mathrm x + \mathrm u + \mathrm w)}$} by removing linearly dependent rows. 
We use the reordered variables \mbox{$\bar y = P y$} as in~\eqref{eq:step3sol} and write~\eqref{eq:mpcprogpartitioned_coupling} as $\bar D \bar y =  d$ with
\renewcommand*{\arraystretch}{.8}
\begin{IEEEeqnarray*}{l}
\arraycolsep=2.2pt
  \left[\bar D, d \right] = \hspace{-1mm}  
  \left[ \arraycolsep=1.8pt \begin{array}{cccc|c}  
   \external{B}_r &  -W_r &     &                                     &  -\external{A}_r x^1\\ 
                 &       &  F  & 0_{(N-1)\mathrm w \times \mathrm x} & 0\end{array}\right] \in \mathbb R^{N\mathrm w \times (\mathrm y+1)},
  \end{IEEEeqnarray*}
\renewcommand*{\arraystretch}{1}%
where $F$ is the block-diagonal matrix $I_{N-1} \otimes [\external{A}_r, \external{B}_r, -W_r]$.
 We obtain the final form $Dy=d$ in~\eqref{eq:mpcfinal-coupling} with $D=\bar D P$.

\subsection{Justification of Algorithm~\ref{alg:ouradmm} and Proposition~\ref{prop:ouradmmconvergence}}
\label{app:proofouradmm}
 We rewrite the consensus constraint in~\eqref{eq:mpcprogsplitting} as
\begin{IEEEeqnarray}{rCl}
  \textcolor{black}{\tfrac{1}{\sqrt{\rho}} E^{\nicefrac12}}\bmat{I;I}y = \textcolor{black}{\tfrac{1}{\sqrt{\rho}} E^{\nicefrac12}}\bmat{{\zeta} ; {\epsilon}},
  \label{eq:scaledconsensus}
\end{IEEEeqnarray}
where $E \sequal \diag(E_{\zeta}, E_{\epsilon})$, $E_{\zeta}\sequal\diag(E_{\zeta_1},\sdots,E_{\zeta_M})$, and $E_{\epsilon} \sequal \diag(E_{\epsilon_1},\sdots,E_{\epsilon_M})$. The individual scaling matrices are \mbox{$E_{\zeta_i} = \beta \rho_i I_{\mathrm y_i}$} and \mbox{$E_{\epsilon_i} = (1-\beta) \rho_i I_{\mathrm y_i}$}. We require $E$ to be positive definite, which ensures that~\eqref{eq:scaledconsensus} is equivalent to~\eqref{eq:mpcprogsplitting_consenus}. Positive definiteness of $E$ is given if and only if $\rho_i >0$ and $\beta \in (0,1)$. When we use an augmented Lagrangian as in~\cite{raghunathan2014optimal}, we obtain
 \begin{IEEEeqnarray}{rCl}
 \IEEEyesnumber\eqlabel{eq:scaledLagrangianandDefs} \IEEEyessubnumber*
  {\mathcal L}_\rho \big( y, \hspace{-.4ex} \bmat{\zeta\\ \epsilon}\hspace{-.8ex}, {\bar \lambda} \big) &=&
      f(y) + g(\zeta,\epsilon) 
      + \tfrac{1}{2} \| \hspace{-.8ex} \bmat{y \\ y} \sminus \bmat{\zeta\\ \epsilon} \sminus \bmat{{\bar \lambda_\zeta}\\  {\bar \lambda_\epsilon}} \hspace{-.8ex}\|_{{E}}^2 \quad\quad\label{eq:LagrangianScaled}\\
f(y) &=& \sum\nolimits_i \left(\tfrac12 y_i^\top \mathcal Q_i y_i + q_i^\top y_i + \mathcal I_{C_i y_i = c_i}(y_i) \right) \quad \quad\\
g(\zeta,\epsilon) &=& \sum\nolimits_i \mathcal I_{\mathcal Y_i}(\zeta_i) + \mathcal I_{D\epsilon = d}(\epsilon), 
 \end{IEEEeqnarray} 
where
$\bar \lambda = \sqrt{\rho} E^{-\nicefrac12} \lambda$, $\|x\|_E = \|E^{\nicefrac12}x \|_2 = \sqrt{x^\top E x}$, and $\mathcal I$ denotes an indicator function. By abusing the notation, we replace~$\bar \lambda$ with~$\lambda$.
 We obtain Algorithm~\ref{alg:ouradmm} by applying standard ADMM~\cite[Eqn.~(3)]{boyd2011distributed} with~\eqref{eq:scaledLagrangianandDefs}. 
For $\bar \proj_{D\epsilon=d}(\cdot)$, we obtain
 \begin{IEEEeqnarray*}{rCl}
{{\bar \proj}_{{D\epsilon=d}}(\cdot)} = \arg \min\nolimits_{\epsilon} \big\{ \tfrac{1}{2}\| {E_{\epsilon}^{\nicefrac12}}(\cdot - \epsilon)\|_2^2   \,\,\,\text{s.t.}\,\,\, D \epsilon = d \big\},
\end{IEEEeqnarray*}
which results in~\eqref{eq:modifiedProjectione} through substituting~$\bar \epsilon = E_{\epsilon}^{\nicefrac12}\epsilon$.
Proposition~\ref{prop:ouradmmconvergence}-$(i)$ follows as we have reduced Algorithm~\ref{alg:ouradmm} to an application of standard ADMM~\cite{boyd2011distributed}. The algorithm converges to a single fixed point according to~\cite[Thm.~2]{raghunathan2014optimal}, which applies to our formulation as shown in Appendix~\ref{app:proofoptimalparam}. Statement $(ii)$ follows from inserting the parameters. 
\ifArxiv
We provide additional details on the algorithm formulation in the ancillary material.
\else
We provide additional details on the algorithm formulation in~\cite{Rey2018}.
\fi
\hfill $\blacksquare$

\subsection{Justification of  Table~\ref{tab:StepComplexities}}
\label{app:StepComplexities} We denote the cost of an operation with $\cost(\cdot)$.
Based on~\eqref{eq:step1solvepenalty}, the largest subsystem in \circled{2.1} has $\mathcal O(\max_i\cost(\mathcal M_i))$ as $\mathcal N_i c_{i}$ can be precomputed. We use $\cost(\mathcal M_i) \sequal \mathcal O(\cost({P_i}) + \cost({C_i}) + \cost(C_i \mathcal P_i {C_i}^{-1}))$. The cost for multiplying~$P_i$ relates to $N$-times applying the inverse of $Q_i+I$ and $R_i+I$, where the first part dominates due to $\mathrm x \hspace{-0.5ex}\geq\hspace{-0.5ex} \mathrm u$. We precompute an $LDL^\top$ factorization and perform a forward-backward substitution as in~\cite[Sec. 3.1]{golub1996} in~$\mathcal O(N\mathrm x_i^2)$. The cost for multiplying $C_i$ is dominated by $N$ multiplications with $A_{ii}$, which results in~$\mathcal O(N\mathrm x_i^2)$. For $\cost(C_i \mathcal P_i {C_i}^{-1})$, we again use $LDL^\top$, where~$L$ is~$2\mathrm x_i$-banded as the multistage structure makes $C_i$ banded. By following~\cite[Sec. 4.3]{golub1996}, we get~$\mathcal O(N\mathrm x_i^2)$. 

The result for \circled{2.2} follows from the composition of $\mathcal Y$ of $\MN$ times $\mathcal X_i, \mathcal U_i$. In the case of box constraints, the projection reduces to an element-wise clipping in $\mathcal O(2\mathrm x_i + 2\mathrm u_i)=\mathcal O(\mathrm x_i)$ for each subsystem and time instance.  With $2MN\mathrm x_i$ parallel threads, \circled{2.2} can also be executed in $\mathcal O(1)$. 

By following~\eqref{eq:step3sol}, we use $\mathcal O(\cost(\bar{\mathcal D}))$ for \circled{2.3} as we precompute $\bar{\mathcal E}d$ and neglect permutations and sums. We use~$\cost(\bar{\mathcal D}) \sequal \mathcal O(\cost({\bar D}) \splus \cost({(\bar D \bar E_{\epsilon} \bar{D}^\top)}^{-1})$ with~$\bar D$ from Appendix~\ref{app:stacking}. $\bar D$~is dominated by $N$ diagonal blocks $\bar D^k\sequal[\external{A}_r, \external{B}_r, -W_r]$. For the largest block, we require $\mathcal O(\max_k(\cost({\bar D^k}) + \cost((\bar{D}^k \bar{E_\epsilon}^k (\bar{D}^k)^\top)^{-1}))$. We neglect $\cost({\bar D^k})$ as $\bar D^k$  is sparse if the system has structure. For $\cost((\bar{D}^k \bar{E_\epsilon}^k (\bar{D}^k)^\top)^{-1}))$, we use an $LDL^\top$ factorization~\cite[Sec. 3.1]{golub1996} in $\mathcal O(\mathrm w^2)$. In the `out-1' case, $\external{A}, \external{B}$ only have one element per column, hence we can reshuffle $\bar D^k \bar{E_\epsilon}^k ({\bar D^k})^\top$ to become block-diagonal. The largest block then has the size $\max_i \mathrm w_i$.

Step~\circled{2.4} decomposes into $2\MN$ operations of size $\mathrm y_i^k$. The longest thread has $\mathcal O(\max_i\{\mathrm u_i + \mathrm w_i +\mathrm x_i\})\sequal \mathcal O(\max_i \mathrm x_i)$. Further, with $2\MN\mathrm x_i$ parallel threads, it even is $\mathcal O(1)$.

\subsection{Proof of Proposition~\ref{prop:optimal param}}
\label{app:proofoptimalparam}
Our result extends the convergence analysis in~\cite{raghunathan2014optimal} by including a scaled consensus constraint. Instead of $y$, $w$, $Q$, $q$, $A$, $b$, $\mathcal Y$ as used in~\cite{raghunathan2014optimal}, we use the notation $\hat y$, $\hat w$, $\hat Q$, $\hat q$, $\hat A$, $\hat b$, $\hat{\mathcal Y}$, where~\eqref{eq:mpcprogsplitting} relates to~\cite[Eqn.~(3)]{raghunathan2014optimal} through
 \begin{IEEEeqnarray}{rCl}
 \IEEEyesnumber\eqlabel{eq:RaghDefinitions} \IEEEyessubnumber*
 \hat y &=& [\bar y; y], \,\,\hat w = [\zeta; \epsilon], \,\,\hat \lambda = [\lambda_\zeta; \lambda_\epsilon]\\
 \hat Q &=& \nicefrac12 \operatorname{diag}(\mathcal Q, \mathcal Q), \text{ where } \mathcal Q = \diag({\mathcal Q}_1,\sdots,{\mathcal Q}_M)\quad \quad \,\,\,\\
 \hat q &=& \nicefrac12 \, [q;q], \text{ where } q = [q_1; \dots; q_M]\\
 \hat A &=& \bmat{I_{\mathrm y}&-I_{\mathrm y}\\0_{N\mathrm x \times \mathrm y}& C}, \text{ where } C = \diag(C_1,\sdots,C_M)\\
 \hat b &=& [0_{\mathrm y\times 1};c], \text{ where } c = [c_1; \dots; c_M]\\
 \hat{\mathcal Y} &=& \{ \hat w = [\zeta; \epsilon] \, | \, \zeta_i \in \mathcal Y_i, i=1,\sdots, M, D\epsilon = d\}.
 \end{IEEEeqnarray}
 In \cite{raghunathan2014optimal}, it is required that $\{\hat y\, | \, \hat A \hat y = \hat b\} \cap \hat{\mathcal Y}\neq \emptyset$, $\hat A$ has full row rank, and $\hat Z^\top \hat Q \hat Z$ is positive definite, where $\hat Z$ contains an orthonormal null space basis for $\hat A$. 
With~\eqref{eq:RaghDefinitions}, $\hat Z = \tfrac{1}{\sqrt{2}}[Z;Z]$, and the assumptions made for the initial problem and Proposition~\ref{prop:optimal param}, these conditions are satisfied. 
 By using a scaling as in~\eqref{eq:scaledconsensus} and following~\cite{raghunathan2014optimal}, we obtain Algorithm~\ref{alg:eqformtoRaghunathan} and the definitions below.
 \begin{algorithm}[H]
\caption{ADMM as in~\cite[Eqn.~(5)]{raghunathan2014optimal}}
\label{alg:eqformtoRaghunathan}
\begin{algorithmic}[1]
\vspace{0.2em}

\Statex \hspace{-1.05em}\textbf{repeat}
\vspace{0.1em}
\,
\begin{minipage}[t]{\linewidth}
\begin{algorithmic}[1]
\algrenewcommand{\alglinenumber}[1]{\hspace{1.4ex}\footnotesize\color{black}\circled{3.1}}
\State $\hat y\, \leftarrow {\hat{\mathcal{M}}}(\hat w + \hat \lambda- {E^{-1}} \hat q) + {\hat{\mathcal{N}}} \hat b$
\vspace{0.4em}
\algrenewcommand{\alglinenumber}[1]{\hspace{1.4ex}\footnotesize\color{black}\circled{3.2}}
\State $\hat w \leftarrow {T(  \textcolor{black}{\hat y -\hat \lambda})}$ 
\vspace{0.4em}
\algrenewcommand{\alglinenumber}[1]{\hspace{1.4ex}\footnotesize\color{black}\circled{3.3}}
\State $  \hat \lambda \,\leftarrow  \hat \lambda  -  \hat y +   \hat w$
\end{algorithmic}
\end{minipage}
\end{algorithmic}
\end{algorithm}%
\vspace{-1.8\baselineskip}
\begin{IEEEeqnarray}{rCl}
\IEEEyesnumber\eqlabel{eq:hatMNandTdefs} \IEEEyessubnumber*
  \hat{\mathcal M} &=& \bmat{I_{\mathrm y} ; I_{\mathrm y}}\diag({\mathcal M}_i) \big[{\diag\big(\tfrac{1}{\rho_i}E_{\zeta_i}\big)}, {\diag\big(\tfrac{1}{\rho_i}E_{\epsilon_i}\big)} \big]\label{eq:hatMdef}\quad\quad\\
  \hat{\mathcal N} &=& \bmat{0_{(\mathrm y + N\mathrm x)\times \mathrm y} &\hspace{-1ex} \diag({\mathcal N}_i)}\\
  T(\cdot) &=& {E^{-\nicefrac12}}\proj\nolimits_{{E^{\nicefrac12}}\hat{\mathcal Y}}({E^{\nicefrac12}}\cdot)
\end{IEEEeqnarray}
In~\eqref{eq:hatMNandTdefs}, we use ${\mathcal M}_i$ and ${\mathcal N}_i$ as in~\eqref{eq:step1solvepenalty}. To be consistent with~\cite{raghunathan2014optimal}, we require $(i)$ that $\hat{\mathcal M}$ is a contraction, and $(ii)$ that $T(\cdot)$ is firmly nonexpansive~\cite{bauschke2017convex}, which is equivalent to the properties noted in~\cite[Eqn.~(8), Lem.~3]{raghunathan2014optimal} by~\cite[Definition~4.1, Proposition~4.2]{bauschke2017convex}.
For $(i)$, it is sufficient to check if all $\mathcal M_i$ are contractions. 
Similar to~\cite[Eqn.~(5)]{raghunathan2014optimal}, we use the null space method in~\cite[Sec. 16.2]{Nocedal2000} to obtain the equivalent form $\mathcal M_i =  Z_i( Z_i^\top(\frac1{\rho_i}  Q_i+I) Z_i)^{-1} Z_i^\top$. It can be shown that $(i)$ is satisfied if $\hat Z_i^\top \hat Q_i \hat Z_i$ is positive definite.
 Requirement $(ii)$ is true by~\cite[Proposition~4.8]{bauschke2017convex} as~$T(\cdot)$ is an orthogonal projection onto a convex set in a  Hilbert space with inner product $\langle x , y \rangle_{E} \sequal x^\top E y$. As $(i)$, $(ii)$ are satisfied, the convergence analysis in~\cite{raghunathan2014optimal} applies. In~\cite[Sec.~V]{raghunathan2014optimal}, the worst-case convergence rate is optimized with
$\rho^\star \sequal \arg \min\nolimits_{\rho} \| \tilde{\mathcal M}\|_2$,
where $\tilde{\mathcal M} \sequal  \hat Z^\top \hat{\mathcal M} \hat Z - \tfrac12 I_{N\mathrm x}$, which becomes $\tilde{\mathcal M} \sequal \diag ( {Z_i^\top { {\mathcal M}_i}  Z_i - \tfrac12 I_{N\mathrm x_i}})$ in our case.  As~$\tilde{\mathcal M}$ is block-diagonal, we choose separate~$\rho_i$ in the same way as~$\rho$ is chosen in~\cite{raghunathan2014optimal}, which results in Proposition~\ref{prop:optimal param}. 
\ifArxiv
We provide a detailed version of this proof in the ancillary material.
\else
We show a more-detailed version of this proof in~\cite{Rey2018}.
\fi
\hfill $\blacksquare$

\subsection{Proof of Proposition~\ref{prop:stproperties}}
\label{app:proofstproperties}
For $(i)$, we show $0 \hspace{-.65ex}<\hspace{-.65ex} \sum_j \hspace{-.3ex} \Gamma_{ij} \hspace{-.65ex}< \hspace{-.65ex}\infty$ for all~$i$. First, we consider~\eqref{eq:initsys} with $u^{k} \sequal \delta^{k}$, \mbox{$x^{0} \sequal 0$}. We obtain $x^{k} \sequal (A)^{k-1} B 1_{\mathrm u \times 1}$ for $k\geq 1$.
Given that $A$ is semi-convergent, $x^{k}$ asymptotically converges to $x_{eq} \sequal Ax_{eq}$, which implies that $\Delta x^{k}$ asymptotically converges to zero. Hence, $\sum_{k=0}^\infty |\Delta x^{k}|^2$ is a sum over a squared-exponential tail, which is finite~\cite[Eq.~(2151), p.1132]{bronstein2012taschenbuch}. Consequently, each~$\Gamma_{ij}$ is finite and $\sum_j \Gamma_{ij} \hspace{-.3ex} < \hspace{-.3ex} \infty$ for all~$i$. 
To show $\sum_j \Gamma_{ij} \hspace{-.3ex}> \hspace{-.3ex} 0$, we use that controllability implies that $[B, AB, A^2B, \sdots, A^{\mathrm x-1}B]$ has full rank, which means that the sequence $\{x^{k}\} \sequal \{A^{k-1} B 1_{\mathrm u \times 1}\}$ spans~$\mathbb R^{\mathrm x}$. The same is true for $\{\Delta x^{k}\}$. Also, controllability implies that $A,B$ do not have a common zero row. Hence, no row in $\Phi^{k}$ is filled with zeros for all times, and we obtain $\sum_j \Gamma_{ij}>0$ for all $i$.

Statement $(ii)$ is clear from $0 < \sum_j \Gamma_{ij} < \infty$ and~\eqref{eq:defsa}.
 
 For $(iii)$, we use the diagonal state and input transformation $T=\operatorname{diag}(t_1, t_2, \dots, t_{\mathrm x + \mathrm u})$ with $t_i \neq 0$. More specifically, we use $[\bar x^{k},\bar u^{k}] = T^{-1}[x^{k},u^{k}]$. For the transformed system, we obtain $\bar \Phi^{k} = T^{-1}\Phi^{k}$ and $\bar \Gamma_{ij} = |t_i| \Gamma_{ij}$. The factors $|t_i|$ then cancel out in~\eqref{eq:defsa}, which makes $s$ invariant to $T$.\hfill $\blacksquare$


\ifCLASSOPTIONcaptionsoff
  \newpage
\fi



\bibliographystyle{IEEEtran}
\bibliography{IEEEabrv,root}

\begin{thebibliography}{10}
\providecommand{\url}[1]{#1}
\csname url@samestyle\endcsname
\providecommand{\newblock}{\relax}
\providecommand{\bibinfo}[2]{#2}
\providecommand{\BIBentrySTDinterwordspacing}{\spaceskip=0pt\relax}
\providecommand{\BIBentryALTinterwordstretchfactor}{4}
\providecommand{\BIBentryALTinterwordspacing}{\spaceskip=\fontdimen2\font plus
\BIBentryALTinterwordstretchfactor\fontdimen3\font minus
  \fontdimen4\font\relax}
\providecommand{\BIBforeignlanguage}[2]{{%
\expandafter\ifx\csname l@#1\endcsname\relax
\typeout{** WARNING: IEEEtran.bst: No hyphenation pattern has been}%
\typeout{** loaded for the language `#1'. Using the pattern for}%
\typeout{** the default language instead.}%
\else
\language=\csname l@#1\endcsname
\fi
#2}}
\providecommand{\BIBdecl}{\relax}
\BIBdecl

\bibitem{maciejowski2002predictive}
J.~M. Maciejowski, \emph{Predictive Control with Constraints}.\hskip 1em plus
  0.5em minus 0.4em\relax Prentice Hall, 2002.

\bibitem{camacho2013model}
E.~F. Camacho and C.~B. Alba, \emph{Model predictive control}.\hskip 1em plus
  0.5em minus 0.4em\relax Springer Science \& Business Media, 2013.

\bibitem{morari1999model}
M.~Morari and J.~H. Lee, ``Model predictive control: past, present and
  future,'' \emph{Computers \& Chemical Engineering}, vol.~23, no. 4-5, pp.
  667--682, 1999.

\bibitem{di2012industry}
S.~Di~Cairano, ``An industry perspective on {MPC} in large volumes
  applications: Potential benefits and open challenges,'' \emph{IFAC
  Proceedings Volumes}, vol.~45, no.~17, pp. 52--59, 2012.

\bibitem{glowinski1975approximation}
R.~Glowinski and A.~Marroco, ``Sur l'approximation, par {\'e}l{\'e}ments finis
  d'ordre un, et la r{\'e}solution, par p{\'e}nalisation-dualit{\'e} d'une
  classe de probl{\`e}mes de dirichlet non lin{\'e}aires,'' \emph{Revue
  fran{\c{c}}aise d'automatique, informatique, recherche op{\'e}rationnelle.
  Analyse num{\'e}rique}, vol.~9, no.~2, pp. 41--76, 1975.

\bibitem{boyd2011distributed}
S.~Boyd, N.~Parikh, E.~Chu, B.~Peleato, and J.~Eckstein, ``Distributed
  optimization and statistical learning via the alternating direction method of
  multipliers,'' \emph{Foundations and Trends in Machine Learning}, vol.~3,
  no.~1, pp. 1--122, 2011.

\bibitem{jerez2014embedded}
J.~L. Jerez, P.~J. Goulart, S.~Richter, G.~A. Constantinides, E.~C. Kerrigan,
  and M.~Morari, ``Embedded online optimization for model predictive control at
  megahertz rates,'' \emph{IEEE Transactions on Automatic Control}, vol.~59,
  no.~12, pp. 3238--3251, 2014.

\bibitem{raghunathan2014optimal}
A.~U. Raghunathan and S.~Di~Cairano, ``Optimal step-size selection in
  alternating direction method of multipliers for convex quadratic programs and
  model predictive control,'' in \emph{Proceedings of Symposium on Mathematical
  Theory of Networks and Systems}, 2014, pp. 807--814.

\bibitem{jerez2013embedded}
J.~L. Jerez, P.~J. Goulart, S.~Richter, G.~A. Constantinides, E.~C. Kerrigan,
  and M.~Morari, ``Embedded predictive control on an {FPGA} using the fast
  gradient method,'' in \emph{European Control Conference}.\hskip 1em plus
  0.5em minus 0.4em\relax IEEE, 2013, pp. 3614--3620.

\bibitem{peyrl2014parallel}
H.~Peyrl, A.~Zanarini, T.~Besselmann, J.~Liu, and M.-A. Bo{\'e}chat, ``Parallel
  implementations of the fast gradient method for high-speed {MPC},''
  \emph{Control Engineering Practice}, vol.~33, pp. 22--34, 2014.

\bibitem{vaya2014decentralized}
M.~G. Vay{\'a}, G.~Andersson, and S.~Boyd, ``Decentralized control of plug-in
  electric vehicles under driving uncertainty,'' in \emph{Innovative Smart Grid
  Technologies Conference Europe}.\hskip 1em plus 0.5em minus 0.4em\relax IEEE,
  2014, pp. 1--6.

\bibitem{kang2015decomposition}
J.~Kang, A.~U. Raghunathan, and S.~Di~Cairano, ``Decomposition via {ADMM} for
  scenario-based model predictive control,'' in \emph{American Control
  Conference}.\hskip 1em plus 0.5em minus 0.4em\relax IEEE, Jul 2015, pp.
  1246--1251.

\bibitem{Rey2017d}
F.~Rey, X.~Zhang, S.~Merkli, V.~Agliati, M.~Kamgarpour, and J.~Lygeros,
  ``Strengthening the group: Aggregated frequency reserve bidding with admm,''
  \emph{IEEE Transactions on Smart Grid}, 2018.

\bibitem{deng2017parallel}
W.~Deng, M.-J. Lai, Z.~Peng, and W.~Yin, ``Parallel multi-block {ADMM} with
  o(1/k) convergence,'' \emph{Journal of Scientific Computing}, vol.~71, no.~2,
  pp. 712--736, 2017.

\bibitem{chen2016direct}
C.~Chen, B.~He, Y.~Ye, and X.~Yuan, ``The direct extension of {ADMM} for
  multi-block convex minimization problems is not necessarily convergent,''
  \emph{Mathematical Programming}, vol. 155, no. 1-2, pp. 57--79, 2016.

\bibitem{mateos2010distributed}
G.~Mateos, J.~A. Bazerque, and G.~B. Giannakis, ``Distributed sparse linear
  regression,'' \emph{IEEE Transactions on Signal Processing}, vol.~58, no.~10,
  pp. 5262--5276, 2010.

\bibitem{chang2015multi}
T.-H. Chang, M.~Hong, and X.~Wang, ``Multi-agent distributed optimization via
  inexact consensus {ADMM},'' \emph{IEEE Transactions on Signal Processing},
  vol.~63, no.~2, pp. 482--497, 2015.

\bibitem{stellato2017osqp}
B.~Stellato, G.~Banjac, P.~Goulart, A.~Bemporad, and S.~Boyd, ``{OSQP}: An
  operator splitting solver for quadratic programs,'' \emph{arXiv preprint
  arXiv:1711.08013}, 2017.

\bibitem{raghunathan2014admm}
A.~U. Raghunathan and S.~Di~Cairano, ``{ADMM} for convex quadratic programs:
  {Q}-linear convergence and infeasibility detection,'' \emph{arXiv preprint
  arXiv:1411.7288}, 2014.

\bibitem{GhadimiOptimalParameter}
E.~Ghadimi, A.~Teixeira, I.~Shames, and M.~Johansson, ``Optimal parameter
  selection for the alternating direction method of multipliers ({ADMM}):
  Quadratic problems,'' \emph{IEEE Transactions on Automatic Control}, vol.~60,
  no.~3, pp. 644--658, Mar 2015.

\bibitem{Rey2016}
F.~Rey, D.~Frick, A.~Domahidi, J.~Jerez, M.~Morari, and J.~Lygeros, ``{ADMM}
  prescaling for model predictive control,'' in \emph{Conference on Decision
  and Control}.\hskip 1em plus 0.5em minus 0.4em\relax IEEE, 2016, pp.
  3662--3667.

\bibitem{stathopoulos2016operator}
G.~Stathopoulos, H.~Shukla, A.~Szucs, Y.~Pu, C.~N. Jones \emph{et~al.},
  ``Operator splitting methods in control,'' \emph{Foundations and Trends in
  Systems and Control}, vol.~3, no.~3, pp. 249--362, 2016.

\bibitem{giselsson2017linear}
P.~Giselsson and S.~Boyd, ``Linear convergence and metric selection for
  douglas-rachford splitting and {ADMM},'' \emph{IEEE Transactions on Automatic
  Control}, vol.~62, no.~2, pp. 532--544, 2017.

\bibitem{wang2010fast}
Y.~Wang and S.~Boyd, ``Fast model predictive control using online
  optimization,'' \emph{IEEE Transactions on Control Systems Technology},
  vol.~18, no.~2, pp. 267--278, 2010.

\bibitem{Rey2017a}
F.~Rey, P.~Hokayem, and J.~Lygeros, ``A tailored {ADMM} approach for power
  coordination in variable speed drives,'' \emph{20th IFAC World Congress},
  vol.~50, no.~1, pp. 7403--7408, 2017.

\bibitem{Rey2017b}
------, ``Ask not what {ADMM} can do for you, ask what you can do for {ADMM} -
  virtual subsystems in {MPC},'' in \emph{Conference on Decision and
  Control}.\hskip 1em plus 0.5em minus 0.4em\relax IEEE, 2017, pp. 4357--4362.

\bibitem{ageev2000approximation}
A.~A. Ageev and M.~I. Sviridenko, ``An approximation algorithm for hypergraph
  max k-cut with given sizes of parts,'' in \emph{European Symposium on
  Algorithms}.\hskip 1em plus 0.5em minus 0.4em\relax Springer, 2000, pp.
  32--41.

\bibitem{frieze1997improved}
A.~Frieze and M.~Jerrum, ``Improved approximation algorithms for maxk-cut and
  max bisection,'' \emph{Algorithmica}, vol.~18, no.~1, pp. 67--81, 1997.

\bibitem{hespanha2004efficient}
J.~P. Hespanha, ``An efficient matlab algorithm for graph partitioning,''
  \emph{University of California}, pp. 1--8, 2004.

\bibitem{parikh2014proximal}
N.~Parikh, S.~Boyd \emph{et~al.}, ``Proximal algorithms,'' \emph{Foundations
  and Trends in Optimization}, vol.~1, no.~3, pp. 127--239, 2014.

\bibitem{banjac2017infeasibility}
G.~Banjac, P.~Goulart, B.~Stellato, and S.~Boyd, ``Infeasibility detection in
  the alternating direction method of multipliers for convex optimization,''
  \emph{Optimization Online}, 2017.

\bibitem{Nocedal2000}
J.~Nocedal and S.~J. Wright, \emph{Numerical optimization}.\hskip 1em plus
  0.5em minus 0.4em\relax Springer Science and Business Media, 1975, vol.~9,
  no.~4.

\bibitem{domahidi2012efficient}
A.~Domahidi, A.~U. Zgraggen, M.~N. Zeilinger, M.~Morari, and C.~N. Jones,
  ``Efficient interior point methods for multistage problems arising in
  receding horizon control,'' in \emph{Conference on Decision and
  Control}.\hskip 1em plus 0.5em minus 0.4em\relax IEEE, 2012, pp. 668--674.

\bibitem{okuyama2014discrete}
Y.~Okuyama, \emph{Discrete control systems}.\hskip 1em plus 0.5em minus
  0.4em\relax Springer, 2014.

\bibitem{alyani2017diagonality}
K.~Alyani, M.~Congedo, and M.~Moakher, ``Diagonality measures of hermitian
  positive-definite matrices with application to the approximate joint
  diagonalization problem,'' \emph{Linear Algebra and its Applications}, vol.
  528, pp. 290--320, 2017.

\bibitem{toscano2013structured}
R.~Toscano, \emph{Structured controllers for uncertain systems}.\hskip 1em plus
  0.5em minus 0.4em\relax Springer, 2013.

\bibitem{oppenheim1999discrete}
A.~V. Oppenheim, \emph{Discrete-time signal processing}.\hskip 1em plus 0.5em
  minus 0.4em\relax Pearson Education India, 1999.

\bibitem{gu2012discrete}
G.~Gu, \emph{Discrete-Time Linear Systems: Theory and Design with
  Applications}.\hskip 1em plus 0.5em minus 0.4em\relax Springer Science \&
  Business Media, 2012.

\bibitem{cantoni2017structured}
M.~Cantoni, F.~Farokhi, E.~Kerrigan, and I.~Shames, ``Structured computation of
  optimal controls for constrained cascade systems,'' \emph{International
  Journal of Control}, pp. 1--10, 2017.

\bibitem{li2005water}
Y.~Li, M.~Cantoni, and E.~Weyer, ``On water-level error propagation in
  controlled irrigation channels,'' in \emph{Conference on Decision and
  Control}.\hskip 1em plus 0.5em minus 0.4em\relax IEEE, 2005, pp. 2101--2106.

\bibitem{soltanian2015decentralized}
L.~Soltanian and M.~Cantoni, ``Decentralized string-stability analysis for
  heterogeneous cascades subject to load-matching requirements,''
  \emph{Multidimensional Systems and Signal Processing}, vol.~26, no.~4, pp.
  985--999, 2015.

\bibitem{labadie2004optimal}
J.~W. Labadie, ``Optimal operation of multireservoir systems: state-of-the-art
  review,'' \emph{Journal of water resources planning and management}, vol.
  130, no.~2, pp. 93--111, 2004.

\bibitem{guo2011hierarchical}
G.~Guo and W.~Yue, ``Hierarchical platoon control with heterogeneous
  information feedback,'' \emph{IET control theory \& applications}, vol.~5,
  no.~15, pp. 1766--1781, 2011.

\bibitem{lofberg2004yalmip}
J.~Lofberg, ``{YALMIP}: A toolbox for modeling and optimization in {MATLAB},''
  in \emph{International Symposium on Computer Aided Control Systems
  Design}.\hskip 1em plus 0.5em minus 0.4em\relax IEEE, 2004, pp. 284--289.

\bibitem{gurobi}
\BIBentryALTinterwordspacing
I.~Gurobi~Optimization, ``Gurobi optimizer reference manual,'' 2016. [Online].
  Available: \url{http://www.gurobi.com}
\BIBentrySTDinterwordspacing

\bibitem{golub1996}
G.~H. Golub and C.~F. Van~Loan, \emph{Matrix Computations}.\hskip 1em plus
  0.5em minus 0.4em\relax Baltimore, MD, USA: Johns Hopkins University Press,
  1996.

\bibitem{bauschke2017convex}
H.~H. Bauschke and P.~L. Combettes, \emph{Convex analysis and monotone operator
  theory in Hilbert spaces}.\hskip 1em plus 0.5em minus 0.4em\relax Springer,
  2017, vol. 2011.

\bibitem{bronstein2012taschenbuch}
I.~N. Bronstein, J.~Hromkovic, B.~Luderer, H.-R. Schwarz, J.~Blath, A.~Schied,
  S.~Dempe, G.~Wanka, and S.~Gottwald, \emph{Taschenbuch der mathematik}.\hskip
  1em plus 0.5em minus 0.4em\relax Springer-Verlag, 2012, vol.~1.

\end{thebibliography}






\begin{IEEEbiography}[{\includegraphics[width=1in,height=1.25in,clip,keepaspectratio]{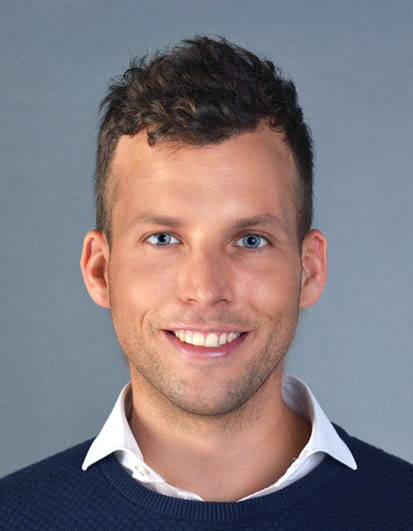}}]{Felix Rey} received his B.Eng.~degree from the University of Applied Sciences Constance (Germany, 2011), and M.Sc.~degree  from the Karlsruhe Institute of Technology (Germany, 2014), both in Electrical Engineering and Information Technology. He is a fellow of the German National Academic Foundation. In 2014, he joined the Automatic Control Laboratory at ETH Zurich as a Ph.D.~candidate. His research interests include model predictive control, as well as distributed and embedded optimization, in particular with ADMM.
\end{IEEEbiography}

 
\begin{IEEEbiography}[{\includegraphics[width=.95in,height=1.25in,clip,keepaspectratio]{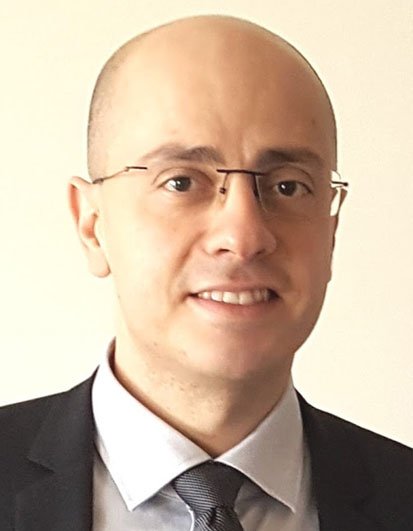}}]{Peter Hokayem}
 received his B.Eng.~degree in Computer and Communications Engineering from the American University of Beirut in 2001, his M.Sc.~degree in Electrical Engineering from the University of New Mexico in 2003, and his Ph.D.~degree in Electrical and Computer Engineering from the University of Illinois at Urbana-Champaign in 2007. Following a Postdoctoral Scholar position at the University of Wuerzburg, he joined the Automatic Control Laboratory at ETH Zurich in 2008 as a (Senior) Postdoctoral Researcher. He moved to ABB Switzerland in 2011 and since then has held several positions, both at the Business Unit and the Corporate Research Center. Currently, he has a consulting role at ABB as a Drive System Expert in Medium Voltage Drives, with a strong focus on novel estimation and control methods for large-scale electromechanical power conversion systems. He was a recipient of the Automatica Paper Prize Award in 2008 and is a Senior Member of the IEEE.
\end{IEEEbiography}

 
\begin{IEEEbiography}[{\includegraphics[width=1in,height=1.25in,clip,keepaspectratio]{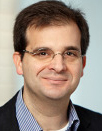}}]{John Lygeros} (M'90--F'11)
  received his B.Eng.~degree in electrical engineering (1990) and his M.Sc.~degree in Systems and Control (1991), both at Imperial College of Science Technology and Medicine, London, UK, and his Ph.D.~degree from the Electrical Engineering and Computer Sciences Department, University of California, Berkeley (1996). He holds the chair of Computation and Control at the Swiss Federal Institute of Technology (ETH) Zurich, Switzerland, where he is currently serving as the Head of the Automatic Control Laboratory. He held a series of research appointments at the National Automated Highway Systems Consortium, Berkeley, the Laboratory for Computer Science, M.I.T., and the Electrical Engineering and Computer Sciences Department at U.C. Berkeley (1996-2000). He was a University Lecturer at the Department of Engineering, University of Cambridge, UK, and a Fellow of Churchill College (2000-2003). He was an Assistant Professor at the Department of Electrical and Computer Engineering, University of Patras, Greece (2003-2006). In 2006 he joined the Automatic Control Laboratory at ETH Zurich, first as an Associate Professor, and since 2010 as Full Professor. His research interests include modeling, analysis, and control of hierarchical, hybrid, and stochastic systems, with applications to biochemical networks, automated highway systems, air traffic management, power grids, and camera networks. He is a Fellow of the IEEE and a Member of the IET and the Technical Chamber of Greece. Since 2013, he is serving as the Treasurer and a Council Member of the International Federation of Automatic Control.
\end{IEEEbiography}

\end{document}